\documentclass[reqno,draft]{amsart}

\theoremstyle{plain}
\newtheorem{lemma}{Lemma}
\newtheorem{theorem}[lemma]{Theorem}

\theoremstyle{remark}
\newtheorem{remark}{Remark}

\def\tt{\theta}

\def\Dd{\Delta}
\def\Gg{\Gamma}
\def\Om{\Omega}

\def\pp{\partial}

\begin{document}

\vskip 0.125in

\title[Primitive Equations]
{Global Well--posedness of the $3D$
 Primitive Equations With Partial Vertical
Turbulence Mixing Heat Diffusion}

\date{October 25, 2010}
\thanks{\textit{ }}

\author[C. Cao]{Chongsheng Cao}
\address[C. Cao]
{Department of Mathematics  \\
Florida International University  \\
University Park  \\
Miami, FL 33199, USA.} \email{caoc@fiu.edu}

\author[E.S. Titi]{Edriss S. Titi}
\address[E.S. Titi]
{Department of Mathematics \\
and  Department of Mechanical and  Aerospace Engineering \\
University of California \\
Irvine, CA  92697-3875, USA. {\bf Also:}
 Department of Computer Science and Applied Mathematics \\
Weizmann Institute of Science  \\
Rehovot 76100, Israel.} \email{etiti@math.uci.edu}
\email{edriss.titi@weizmann.ac.il}

\begin{abstract}

The three--dimensional incompressible viscous
Boussinesq equations, under the assumption of
hydrostatic balance, govern the large scale
dynamics of atmospheric and oceanic motion, and are
commonly called the primitive equations. To
overcome the turbulence mixing a partial vertical
diffusion is usually added to the temperature
advection (or density stratification) equation. In
this paper we prove the  global regularity of
strong solutions to this model in a
three-dimensional infinite horizontal channel,
subject to periodic boundary conditions in the
horizontal directions, and with  no-penetration and
stress-free boundary conditions on the solid, top
and bottom, boundaries. Specifically, we show that
short time strong solutions to the above problem
exist globally in time, and that they depend
continuously on the initial data.
\end{abstract}

\maketitle



\section{Introduction}   \label{S-1}

The partial differential equation model that
describes convective flow in ocean dynamics is known to be the Boussinesq
equations, which are the Navier--Stokes equations (NSE) of incompressible flows with rotation
coupled to the heat (or density stratification) and salinity  transport equations.
The questions of the global
well--posedness of the $3D$ Navier--Stokes equations are considered
to be among the most challenging mathematical problems. In the
context of the atmosphere and the ocean circulation dynamics
geophysicists take advantage of the shallowness of the oceans and
the atmosphere to simplify the Boussinesq equations by modeling the
vertical motion with the hydrostatic balance. This leads to the well-known
primitive equations for ocean and atmosphere dynamics (see, e.g.,
\cite{LTW92}, \cite{LTW92A}, \cite{PJ87}, \cite{Richardson},
\cite{SALMON}, \cite{TZ04}, \cite{VG06} and references therein).
A vertical heat diffusivity is usually added as a
leading order approximation to the effect of micro-scale turbulence mixing
(cf., e.g., \cite{GA84}, \cite{GC93}, \cite{LTW92}, \cite{Richardson}).
As a result one arrives
to the following dimensionless $3D$
 variant of the primitive equations
(Boussinesq equations):
\begin{eqnarray}
&&\hskip-.8in \frac{\pp v}{\pp t} + (v\cdot {\nabla}_{H}) v + w
\frac{\pp v}{\pp z} + f_0 \vec{k} \times v +
{\nabla}_{H} p  + L_1v = 0  \label{EQ-1}  \\
&&\hskip-.8in
\pp_z p  +  T =0,    \label{EQ-2}  \\
&&\hskip-.8in
{\nabla}_{H} \cdot v +\pp_z w =0,   \label{EQ-3} \\
&&\hskip-.8in \frac{\pp T}{\pp t}  + v \cdot {\nabla}_{H} T + w
\frac{\pp T}{\pp z}  +  L_2 T =  Q, \label{EQ-4}
\end{eqnarray}
where the horizontal velocity vector field $v=(v_1, v_2)$,  the velocity vector
field $(v_1, v_2, w)$, the temperature $T$ and the
pressure $p$ are the unknowns. $f_0$ is the Coriolis parameter,
$Q$ is a given heat
source. For simplicity, we drop the coupling with the salinity equation, which is
an advection diffusion equation, but the results reported here will be equally valid
with the addition of  the coupling with the salinity. Moreover, we also assume for simplicity
that $Q$ is time independent.
The viscosity and the heat vertical diffusion operators
$L_1$ and $L_2$, respectively, are given by
\begin{eqnarray}
&&\hskip-.8in
 L_1 = -\frac{1}{ R_1} {\Dd}_{H} -  \frac{1}{  R_2 } \;
\frac{\pp^2}{\pp z^2},   \label{L-1} \\
&&\hskip-.8in
 L_2 =  -  \frac{1}{  R_3} \;
\frac{\pp^2}{\pp z^2},  \label{L-2}
\end{eqnarray}
where $R_1, R_2$ are positive constants representing the horizontal
and vertical dimensionless  Reynolds numbers, respectively, and  $ R_3$ is positive
constant which stands for the vertical dimensionless  eddy heat
diffusivity turbulence mixing coefficient (cf., e.g., \cite{GA84}, \cite{GC93}). We set ${\nabla}_{H}  =
(\pp_x, \pp_y)$ to be the horizontal gradient operator and
${\Dd}_{H} = \pp_x^2 +\pp_y^2$ to be the horizontal Laplacian. We
denote by
\begin{eqnarray}
&&\hskip-.8in
\Gg_u = \{ (x,y,0)  \in  \mathbb{R}^3 \},  \\
&&\hskip-.8in
\Gg_b = \{ (x,y, -h) \in  \mathbb{R}^3 \},
\end{eqnarray}
the upper and lower solid boundaries, respectively.
We equip system (\ref{EQ-1})--(\ref{EQ-4}), on the physical top and bottom boundaries,
with the following  no--normal flow and stress free boundary conditions for the flow velocity
 vector field $(v, w)$, namely,
\begin{eqnarray}
&&\hskip-.8in \mbox{on } \Gg_u: \frac{\pp v }{\pp z} = 0,  \; w=0,
\label{B-1}\\
&&\hskip-.8in \mbox{on } \Gg_b: \frac{\pp v }{\pp z} = 0, \; w=0,   \label{B-2}
\end{eqnarray}
and for simplicity, we set the Dirichlet boundary condition for $T$:
\begin{eqnarray}
&&\hskip-.8in  \left. T\right|_{z=0} =0, \quad \left. T\right|_{z=-h} =1. \label{B-333}
\end{eqnarray}
Horizontally, we set $(v, w)$ and $T$ to satisfy periodic boundary
conditions:
\begin{eqnarray}
&&\hskip-.8in v(x+1,y,z)=v(x,y+1,z)=v(x,y,z);  \label{B-31} \\
&&\hskip-.8in w(x+1,y,z)=w(x,y+1,z)=w(x,y,z);  \label{B-32} \\
&&\hskip-.8in T(x+1,y,z)=T(x,y+1,z)=T(x,y,z).  \label{B-33}
\end{eqnarray}
We will denote by
\[
M=(0, 1)^2
\qquad
\text{and} \quad \Om= M \times (-h, 0).
\]
In addition, we supply the system with the initial condition:
\begin{eqnarray}
&&\hskip-.8in
v(x,y,z,0) = v_0 (x,y,z),  \label{INIT-1}\\
&&\hskip-.8in T(x,y,z,0) = T_0 (x,y,z).  \label{INIT-2}
\end{eqnarray}
System (\ref{EQ-1})--(\ref{INIT-2}) is a modified form of the rotational Rayleigh--B\'{e}nard convection
problem taking into consideration the geophysical situation of the shallowness of oceans and atmosphere.
The original three-dimensional Rayleigh--B\'{e}nard convection model (which is identical, in the absence of
heat diffusion, to the Boussinesq model of stratified fluid)
 has been a subject to study for
many years, numerically, experimentally and analytically (see, e.g.,
 \cite{BH01}, \cite{CC89}, \cite{ES94}, \cite{GE98}, \cite{NSSD}, \cite{PJ87},
 \cite{RL16}, \cite{TZ04}, and references therein). However, the question of global regularity
 is still open and is as challenging as the $3D$ NSE.
Recently, the authors of \cite{CH05} and
\cite{HL05} have shown the global well-posedness to
the $2D$ Boussinesq equations without diffusivity
in the heat transport equation (see also,
\cite{DP08}, for recent improvement). In
\cite{CT05} it is observed that thanks to
hydrostatic balance (\ref{EQ-2}) the unknown
pressure is essentially a function of the
two--dimensional horizontal variables. We take
advantage of this observation in
 \cite{CT05} to establish the $L^6$
estimates for the velocity vector field which allows us to prove the global
well--posedness of the $3D$ primitive equations under the geophysical boundary conditions.
 In \cite{KZ07} the authors take advantage of this
observation as well, and proved the global
well--posedness of the $3D$ primitive equations
with the Dirichlet boundary conditions by dealing
directly with the ``pressure" which is a function
of two variables. In this paper we study system
(\ref{EQ-1})--(\ref{INIT-2}), exploring again the
hydrostatic balance which  leads to an unknown
``pressure" that is a function of only two
variables, and use the techniques and ideas
introduced in \cite{CT05}, \cite{CH05} and
\cite{HL05},   to show in section \ref{S-3} that
strong solutions exist globally in time provided
they exist for a short interval of time.
Furthermore, we show in section \ref{S-4} the
uniqueness and continuous dependence on initial
data of these strong solutions. The short time
existence of strong solutions to this  model will
be reported in a forthcoming paper.

This paper is organized as follows. In section \ref{S-2}, we
reformulate system (\ref{EQ-1})--(\ref{INIT-2}) and introduce
our notations and recall some well-known inequalities.  Section
\ref{S-3} is the main section in which we establish the required
 estimates for proving the global existence in time
 for any initial data. In section \ref{S-4} we prove
the uniqueness of the solutions and their
continuous dependence on initial data.


\section{Functional setting and Formulation}    \label{S-2}

\subsection{Equivalent Formulation}

We denote by
\begin{equation}
\overline{\phi} (x, y)= \frac{1}{h} \int_{-h}^0 \phi(x,y, z)
dz, \qquad \forall \; (x,y) \in M; \label{VBAR}
\end{equation}
and denote the fluctuation by
\begin{equation}
\widetilde{\phi} = \phi - \overline{\phi}.    \label{V--T}
\end{equation}
Notice that
\begin{eqnarray}
&&\hskip-.68in \overline{\widetilde{\phi}} = 0. \label{ZERO}
\end{eqnarray}
Similar to \cite{CT05}, by integrating (\ref{EQ-2}) and (\ref{EQ-3}) vertically, we get
\begin{equation}
w(x,y,z,t) =  - \int_{-h}^z \nabla_H \cdot v(x,y, \xi,t) d\xi,
\label{DIV-1}
\end{equation}
and
\begin{equation}
p(x,y,z,t) = - \int_{-h}^z T(x,y,\xi,t) d\xi + p_s(x,y,t), \label{PPP}
\end{equation}
where $p_s$ is the pressure on the bottom $z=-h$.
Essentially, $p_s(x,y,t)$ is the unknown pressure,
and we observe, as before, that it is a function of
two spatial variables $(x,y).$ As we mentioned in
the introduction we explore this property as in
\cite{CT05} (see also \cite{KZ07}) to prove our
global regularity result.

Replacing $T$ by $T + \frac{z}{h}$, we have the
following equivalent formulation for system
(\ref{EQ-1})--(\ref{INIT-2}):
\begin{eqnarray}
&&\hskip-.68in \frac{\pp v}{\pp t} + L_1 v+ (v\cdot {\nabla}_{H}) v
- \left( \int_{-h}^z {\nabla}_{H} \cdot v(x,y, \xi,t) d\xi
\right)  \frac{\pp v}{\pp  z}  \nonumber  \\
&&\hskip-.5in + f_0 \vec{k} \times v +  {\nabla}_{H} p_s(x,y,t) - {\nabla}_{H} \int_{-h}^z
T(x,y,\xi,t) d\xi = 0,
\label{EQV}   \\
&&\hskip-.68in {\nabla}_{H} \cdot \overline{v} = 0,    \label{EQ22}   \\
&&\hskip-.68in
 \frac{\pp T}{\pp t}  + L_2 T   + v  \cdot {\nabla}_{H} T
- \left( \int_{-h}^z {\nabla}_{H} \cdot v (x,y, \xi,t)
 d\xi
\right) \left( \frac{\pp T}{\pp z}+ \frac{1}{h}\right)   = Q,
 \label{EQ5}  \\
&&\hskip-.68in \left. \frac{\pp v }{\pp z} \right|_{z=0} = \left.
\frac{\pp v }{\pp z} \right|_{z=-h} = 0, \quad v(x+1,y,z)=v(x,y+1,z)=v(x,y,z),
\label{EQ6} \\
&&\hskip-.68in \left.   T \right|_{z=0}=  \left. T \right|_{z=-h}=
0, \quad T(x+1,y,z)=T(x,y+1,z)=T(x,y,z),
 \label{EQ7} \\
&&\hskip-.68in v (x,y,z,0) = v_0 (x,y,z),
\label{EQ8}   \\
&&\hskip-.68in T(x,y,z,0) = T_0 (x,y,z)-\frac{z}{h}.
 \label{EQ9}
\end{eqnarray}
In addition,   $\overline{v}$ and $\widetilde{v}$  satisfy the following coupled system of equations:
\begin{eqnarray}
&&\hskip-.68in \frac{\pp \overline{v}}{\pp t} - \frac{1}{R_1} {\Dd}_{H}
\overline{v} + (\overline{v} \cdot {\nabla}_{H} ) \overline{v} +
\overline{ \left[ (\widetilde{v} \cdot {\nabla}_{H}) \widetilde{v}
 + ({\nabla}_{H} \cdot \widetilde{v}) \; \widetilde{v}\right]}
   \nonumber  \\
&&\hskip-.28in + f_0 \vec{k} \times \overline{v} + {\nabla}_{H} \left[ p_s(x,y,t) - \frac{1}{h}
\int_{-h}^0 \int_{-h}^z  \, T (x,y,\xi,t) \;  d\xi \; dz \right]
= 0,  \label{EQ1}  \\
&&\hskip-.68in {\nabla}_{H} \cdot \overline{v} = 0,    \label{EQ2}
\\
&&\hskip-.68in \frac{\pp \widetilde{v}}{\pp t} + L_1 \widetilde{v} +
(\widetilde{v} \cdot {\nabla}_{H}) \widetilde{v} - \left(
\int_{-h}^z {\nabla}_{H} \cdot \widetilde{v}(x,y, \xi,t) d\xi
\right) \frac{\pp \widetilde{v}}{\pp  z} +(\widetilde{v} \cdot
{\nabla}_{H} )
\overline{v}+ (\overline{v} \cdot {\nabla}_{H}) \widetilde{v}  \nonumber    \\
&&\hskip-.58in   - \overline{ \left[(\widetilde{v} \cdot
{\nabla}_{H}) \widetilde{v}  + ({\nabla}_{H} \cdot \widetilde{v}) \;
\widetilde{v}\right]} + f_0 \vec{k} \times \widetilde{v}  -  {\nabla}_{H} \left( \int_{-h}^z
T(x,y,\xi,t) d\xi -\frac{1}{h} \int_{-h}^0 \int_{-h}^z T(x,y,\xi,t)
d\xi dz \right)  =0.   \label{EQ4}
\end{eqnarray}

\vskip0.1in

\subsection{Functional spaces and inequalities}

Let us denote by $L^q(\Om), L^q(M)$ and $W^{m, q}(\Om), W^{m,
q}(M)$, and $H^m(\Om) =: W^{m, 2}(\Om), H^m(M) =: W^{m, 2}(M)$,
the usual $L^q-$Lebesgue and Sobolev spaces, respectively
(\cite{AR75}). We denote by
\begin{equation}
\| \phi\|_q = \left\{ \begin{array}{ll} \left(  \int_{\Om}  |\phi
(x,y,z)|^q \; dxdydz
\right)^{\frac{1}{q}},   \qquad & \mbox{ for every $\phi \in L^q(\Om)$} \\
\left(  \int_{M}  |\phi (x,y)|^q \; dxdy \right)^{\frac{1}{q}},  &
\mbox{ for every $\phi \in L^q(M)$}.
\end{array} \right.
 \label{L2}
\end{equation}

For convenience, we recall the following Sobolev and
Ladyzhenskaya type inequalities in $(\mathbb{R}/\mathbb{Z})^2$ and in
$\Om$ (see, e.g., \cite{AR75}, \cite{CF88}, \cite{GA94},
\cite{LADY})
\begin{eqnarray}
&&\hskip-.68in
\| \phi \|_{L^4(M)} \leq C_0 \| \phi \|_{L^2}^{1/2} \| \phi \|_{H^1(M)}^{1/2},
\qquad  \forall  \phi \in H^1(M),  \label{SI-1}\\
&&\hskip-.68in \| \phi\|_{L^8(M)} \leq C_0 \| \phi
\|_{L^6(M)}^{3/4} \| \phi \|_{H^1(M)}^{1/4}, \qquad  \forall  \phi \in H^1(M),
\label{SI-2}   \\
&&\hskip-.68in \|{\nabla}_{H} \phi \|_{L^4(M)} \leq C_0 \| \phi
\|_{\infty}^{1/2} \| \phi \|_{H^2(M)}^{1/2},  \qquad  \forall  \phi
\in H^2(M),
\label{SI-11}\\
&&\hskip-.68in \|{\nabla}_{H} \phi \|_{L^4(M)} \leq C_0 \| \phi
\|_{L^2(M)}^{1/2} \| {\nabla}_{H} \phi \|_{\infty}^{1/2} +\| \phi
\|_{L^2(M)}, \qquad  \forall \; \phi \; \mbox{such that } \; {\nabla}_{H} \phi \in L^{\infty}(M),
\label{SI-111}
\end{eqnarray}
and
\begin{eqnarray}
&&\hskip-.68in
\| \psi \|_{L^3(\Om)} \leq C_0 \| \psi \|_{L^2(\Om)}^{1/2} \| \psi \|_{H^1(\Om)}^{1/2},
\label{SI1}\\
&&\hskip-.68in \| \psi \|_{L^6(\Om)} \leq C_0 \| \psi
\|_{H^1(\Om)}, \label{SI2}
\end{eqnarray}
for every $\psi\in H^1(\Om).$  Here $C_0$ is a positive scale invariant constant.
Also, we recall the following version of  Helmholtz-–Weyl decomposition Theorem
(cf. for example, \cite{BM02},  \cite{GA94}, \cite{YV63})
\begin{eqnarray}
&& \|{\nabla}_{H} \phi \|_{W^{m,q}(M)} \leq C \left(
\|{\nabla}_{H} \cdot \phi \|_{W^{m,q}(M)} + \|{\nabla}_{H}^{\perp}
\cdot \phi \|_{W^{m,q}(M)} \right), \label{DIV-CUR}
\end{eqnarray}
for every $\vec{\phi}\in \left( W^{m, q}(M)\right)^2.$ Moreover, we
recall the following Brezis--Gallouet or, Brezis--Wainger inequality (see, e.g.,
\cite{BM02}, \cite{BG80}, \cite{BW80}, \cite{EN89})
\begin{eqnarray}
&& \|\phi \|_{L^{\infty}(M)} \leq C  \|\phi \|_{H^1(M)}\; \left(1+
\log^+ \|\phi \|_{H^2(M)} \right)^{1/2}, \label{BW-1}
\end{eqnarray}
for every $\phi \in H^2(M),$ where $\log^+ r = \log r,$ when $r \geq
1$, and $\log^+ r = 0,$ when $r \leq 1$. Also, we recall the following inequality (see, e.g.,
\cite{BKM84} and \cite{KATO83})
\begin{eqnarray}
&& \|\nabla_H \phi \|_{L^{\infty}(M)} \leq C
(\|\nabla_H \cdot \phi \|_{L^{\infty}(M)} + \|\nabla_H \times \phi \|_{L^{\infty}(M)}) \;
\left(1+ \log^+ \|\nabla_H \phi \|_{H^2(M)} \right), \label{BW-2}
\end{eqnarray}
for every $\nabla_H \phi \in H^2(M)$.
 Moreover, by (\ref{SI-1}) we get
\begin{eqnarray}
&&\hskip-.68in \| \phi \|_{L^{4q}(M)}^{4q} = \|\;
|\phi|^q\;\|_{L^{4}(M)}^{4}
\leq C \|\; |\phi|^q \;\|_{L^2(M)}^{2} \| \;|\phi|^q \;\|_{H^1(M)}^{2}  \nonumber  \\
&&\hskip-.68in \leq C_q  \| \phi \|_{2q}^{2q} \; \left( \int_M
|\phi|^{2q-2} \left| {\nabla}_{H} \phi \right|^2 \; dxdy  \right) +
\| \phi \|_{2q}^{4q},   \label{TWE}
\end{eqnarray}
for every $\phi$ satisfying $\int_M |\phi|^{2q-2} \left| {\nabla}_{H}
\phi \right|^2 \; dxdy < \infty$ and $q\geq 1$. Also, we recall the integral
version of Minkowsky inequality for the $L^p$ spaces, $p\geq 1$. Let
$\Om_1 \subset \mathbb{R}^{m_1}$ and
 $\Om_2 \subset \mathbb{R}^{m_2}$ be two measurable sets, where
$m_1$ and $m_2$ are two positive integers. Suppose that
$f(\xi,\eta)$ is a measurable function over $\Om_1 \times \Om_2$. Then,
\begin{equation}
\hskip0.35in \left[ { \int_{\Om_1} \left( \int_{\Om_2}
|f(\xi,\eta)| d\eta \right)^p d\xi } \right]^{1/p} \leq
\int_{\Om_2} \left( \int_{\Om_1} |f(\xi,\eta)|^p d\xi
\right)^{1/p} d\eta. \label{MKY}
\end{equation}
Finally, we recall the following inequality from Proposition 2.2 in {\bf \cite{CT03}}
\begin{eqnarray}
&&\hskip-.68in
\left| \int_M \left( \int_{-h}^0 \psi_1 (x, y, z)\; dz
\right)\; \left(\int_{-h}^0
\psi_2 (x, y, z) \, \psi_3 (x, y, z)\; dz \right)\; dxdy   \right| \nonumber   \\
&&\hskip-.68in \leq C \| \psi_1 \|_2^{1/2} \| {\nabla}_{H} \psi_1
\|_2^{1/2}\| \psi_2 \|_2^{1/2}\|{\nabla}_{H} \psi_2 \|_2^{1/2}\|
\psi_3 \|_2 +\| \psi_1 \|_2\| \psi_2 \|_2 \| \psi_3 \|_2,
\label{MAIN-1}
\end{eqnarray}
for every $\psi_1, \psi_2 \in H^1(\Om)$ and $\psi_3 \in L^2(\Om),$
and
\begin{eqnarray}
&&\hskip-.68in
\left| \int_M \left( \int_{-h}^0 \psi_1 (x, y, z)\; dz
\right) \, \left( \int_{-h}^0
|{\nabla}_{H} \psi_2(x, y, z)| \, \psi_3(x, y, z)\; dz \right)\; dxdy \right| \nonumber   \\
&&\hskip-.68in  \leq C \| \psi_1 \|_2^{1/2} \| {\nabla}_{H} \psi_1
\|_2^{1/2}\| \psi_2 \|_{\infty}^{1/2}\|{\nabla}_{H}{\nabla}_{H} \psi_2
\|_2^{1/2}\| \psi_3 \|_2 +\| \psi_1 \|_2\| \psi_2 \|_2 \| \psi_3
\|_2, \label{MAIN-2}
\end{eqnarray}
for every $\psi_1 \in H^1(\Om)$, ${\nabla}_{H} \psi_2 \in H^1(\Om)$
and $\psi_3 \in L^2(\Om).$

\section{Global existence of strong solutions} \label{S-3}

In the previous section we have reformulated system
(\ref{EQ-1})--(\ref{INIT-2}) to be equivalent to
system (\ref{EQV})--(\ref{EQ9}). In this section we
will show that strong solutions to system
(\ref{EQV})--(\ref{EQ9})  exist globally in time
provided they exist in short time intervals.

\begin{theorem} \label{T-MAIN}
Let  $Q \in H^2(\Om), v_0 \in H^4(\Om), T_0\in
H^2(\Om)$ and $\mathcal{T} > 0.$ Suppose that there
exists a strong solution $(v(t), T(t))$ of system
{\em (\ref{EQV})--(\ref{EQ9})} on $[0,
\mathcal{T}]$ corresponding to the initial data
$(v_0, T_0)$ such that
\begin{eqnarray*}
&& {\Dd}_{H} v_z, \;\; \nabla_H T \in C([0,\mathcal{T}], H^1(\Om)), \\
&& v_{zz}, {\Dd}_{H} \nabla_H v_z, \;\;
{\nabla}_{H} T_z  \in L^2 ([0,\mathcal{T}], H^1(\Om)).
\end{eqnarray*}
Then this strong solution $(v(t),  T(t))$ exists
globally in time.
\end{theorem}

\begin{remark}
Notice that one can recover the pressure $p_s$ from system (\ref{EQ1})--(\ref{EQ2}) in the same
way as in $2D$ NSE (see, e.g., \cite{CF88}, \cite{SO01A}, \cite{TT84}).
\end{remark}

\begin{proof}

Let $[0,\mathcal{T}_*)$ be the maximal interval of
existence of a strong solution $(v(t), T(t))$. In
order to  establish the global existence, we need
to show  That $\mathcal{T}_* = \infty.$ If
$\mathcal{T}_* < \infty$ we will show  $\|{\Dd}_{H}
v_z(t)\|_{H^1(\Om)}, \;\; \|\nabla_H T
(t)\|_{H^1(\Om)},$ $ \int_0^t\|{\Dd}_{H} \nabla_H
v_z(s)\|_{H^1(\Om)}^2 \; ds $,  $ \int_0^t
\|\nabla_H T_z (s)\|_{H^1(\Om)}^2 \; ds$, and $
\int_0^t \| v_{zz} (s)\|_{H^1(\Om)}^2 \; ds$ are
all bounded uniformly in time, for $t \in
[0,\mathcal{T}_{*})$. As a result the interval $[0,
\mathcal{T}_*)$ can not be a maximal interval of
existence, and consequently the strong solution
$(v(t), T(t))$ exists globally in time.

Therefore, we focus our discussion below on the interval $[0,\mathcal{T}_*).$

\subsection{$\|v\|_2^2+\|T\|_{2}^2$ estimates}

By taking the inner product of equation  (\ref{EQ5}) with $T$, in
$L^2(\Om)$,  we get
\begin{eqnarray}
&&\hskip-.68in \frac{1}{2} \frac{d \|T\|_{2}^{2}}{dt}
+ \frac{1}{R_3}\;\|T_z\|_2^2    \nonumber  \\
&&\hskip-.65in = \int_{\Om} Q T  \; dxdydz -\int_{\Om} \left[ v
\cdot {\nabla}_{H} T - \left( \int_{-h}^z {\nabla}_{H} \cdot v(x,y,
\xi,t) d\xi \right) \left(\frac{\pp T}{\pp z}+\frac{1}{h} \right)
\right]\; T \; dxdydz.   \label{EST-1}
\end{eqnarray}
Integrating by parts and using the boundary
condition (\ref{EQ7}), we get
\begin{eqnarray}
&&\hskip-.065in   -\int_{\Om} \left( v \cdot {\nabla}_{H} T - \left(
\int_{-h}^z {\nabla}_{H} \cdot v(x,y, \xi,t) d\xi \right) \frac{\pp
T}{\pp z}\right) T  \; dxdydz =0.   \label{EST-2}
\end{eqnarray}
As a result of the above we conclude that
\begin{eqnarray}
&&\hskip-.68in \frac{1}{2} \frac{d \|T\|_{2}^{2}}{dt} +
\frac{1}{R_3}\;\|T_z\|_2^2
\nonumber  \\
&&\hskip-.65in =\int_{\Om} \left[ Q  - \frac{1}{h} \left(
\int_{-h}^z {\nabla}_{H} \cdot v(x,y, \xi,t) d\xi \right) \right] \;
T \; dxdydz \leq \|Q\|_{2} \; \|T\|_{2} + \|{\nabla}_{H} v\|_2
\|T\|_2. \label{T-INT}
\end{eqnarray}
Moreover, by taking the inner product of equation
(\ref{EQV}) with $v$, in $L^2(\Om)$,  we reach
\begin{eqnarray}
&&\hskip-.28in \frac{1}{2} \frac{d \|v\|_2^2}{dt} + \frac{1}{R_1}
\|{\nabla}_{H} v\|_2^2
+ \frac{1}{R_2}\|v_z\|_2^2   \nonumber  \\
&&\hskip-.265in =  -\int_{\Om}  \left[  (v \cdot {\nabla}_{H}) v -
\left( \int_{-h}^z {\nabla}_{H} \cdot v(x,y, \xi,t) d\xi \right)
\frac{\pp v}{\pp  z} \right] \cdot  v \;dxdydz
 \nonumber   \\
 &&\hskip-.1065in
- \int_{\Om}  \left( f_0 \vec{k} \times v \right) \cdot v \;dxdydz
- \int_{\Om}  \left(
{\nabla}_{H} p_s - {\nabla}_{H}  \left( \int_{-h}^z T(x,y,\xi,t)
d\xi \right) \right) \cdot v \;dxdydz.  \label{EST-3}
\end{eqnarray}
First, we notice that
\begin{eqnarray}
&&\hskip-.065in \left( f_0 \vec{k} \times v\right) \cdot v =0.   \label{DT-1}
\end{eqnarray}
Next, by integration by parts and using the
boundary conditions (\ref{EQ6}), in particular, the
horizontal periodic boundary conditions, we get
\begin{eqnarray}
&&\hskip-.065in \int_{\Om}  \left[ (v \cdot {\nabla}_{H}) v - \left(
\int_{-h}^z {\nabla}_{H} \cdot v(x,y,\xi,t) d\xi \right) \frac{\pp
v}{\pp z} \right] \cdot v \;dxdydz =0.   \label{DT-11}
\end{eqnarray}
Thanks to (\ref{EQ2}) and, again, the horizontal
periodic boundary conditions, we also have
\begin{eqnarray}
&&\hskip-.065in \int_{\Om} {\nabla}_{H} p_s(x, y, t) \cdot v(x, y,
z, t) \;dxdydz = h \int_M {\nabla}_{H} p_s \cdot \overline{v} \;dxdy
= -h \int_{\Om} p_s ({\nabla}_{H} \cdot \overline{v}) \;dxdy =0.
\label{DT-44}
\end{eqnarray}
By integration by parts, the periodic boundary conditions (\ref{EQ6}), and
Cauchy--Schwarz inequality, we obtain
\begin{eqnarray}
&&\hskip-.065in
\left| \int_{\Om} {\nabla}_{H}  \left( \int_{-h}^z T(x,y,\xi,t)
d\xi \right)  \cdot v \;dxdydz \right| \leq h \, \|T\|_2 \; \|\nabla_H v\|_2.
\label{DT-4}
\end{eqnarray}
Thus, by (\ref{T-INT})--(\ref{DT-4}) we have
\begin{eqnarray}
&&\hskip-.68in \frac{1}{2} \frac{d (\|v \|_2^2+\|T \|_2^2) }{dt} +
\frac{1}{R_1} \|{\nabla}_{H} v\|_2^2
+ \frac{1}{R_2}\|v_z\|_2^2+ \frac{1}{R_3}\|T_z\|_2^2    \nonumber   \\
&&\hskip-.68in  \leq  \|Q\|_2 \|T\|_2+ (1+h) \|T\|_2\; \|
{\nabla}_{H} v \|_2.   \label{EST-4}
\end{eqnarray}
 By Cauchy--Schwarz inequality, we obtain
\begin{eqnarray}
&&\hskip-.68in    \frac{d (\|v \|_2^2+\|T \|_2^2) }{dt}+
\frac{1}{R_1} \|{\nabla}_{H} v\|_2^2
+ \frac{1}{R_2}\|v_z\|_2^2  + \frac{1}{R_3}\|T_z\|_2^2   \\
&&\hskip-.68in \leq  \|Q\|_{2}^2 + (1+R_1)(1+h)^2 \|T\|_2^2.
\label{EST-44}
\end{eqnarray}
Thanks to Gronwall's inequality we get
\begin{eqnarray}
&&\hskip-.68in \| v(t)\|_2^2 +\|T(t) \|_2^2 \leq C \left( \|v_0\|_2^2+
\|T_0\|_2^2 \right) e^{(1+R_1)(1+h)^2 \, t} + C \|Q\|^2_2; \label{V-2}
 \end{eqnarray}
and
\begin{eqnarray}
&&\hskip-.68in \int_0^t \left[ \frac{1}{R_1} \| {\nabla}_{H} v(s)
\|_2^2 + \frac{1}{R_2} \|  v_z (s)\|_2^2 + \frac{1}{R_3}\|T_z(s)\|_2^2
\right]\; ds
 \leq C \left[
  \left( \|v_0\|_2^2+ \|T_0\|_2^2\right)e^{(1+R_1)(1+h)^2 \, t} +
\|Q\|^2_2\; t \right]. \label{VEE}
\end{eqnarray}
Therefore, for every  $t \in [0,\mathcal{T}_*),$ we have
\begin{eqnarray}
&&\hskip-.168in  \| v (t)\|_2^2 +\|T (t)\|_2^2+ \int_0^t \left[ \|
{\nabla}_{H} v (s)\|_2^2  +
 \|  v_z (s)\|_2^2 +\|  T_z (s)\|_2^2  \right]\; ds   \leq K_1,   \label{K-1}
\end{eqnarray}
where
\begin{eqnarray}
&&\hskip-.168in K_1 = C \left[ \left( \|v_0\|_2^2+ \|T_0\|_2^2 \right)
e^{(1+R_1)(1+h)^2 \, t}+ \|Q\|^2_2\; t \right]. \label{K1}
\end{eqnarray}

\subsection{$\|T\|_{\infty}$ estimates}

We follow here the idea of Stampaccia  for proving the Maximum Principle.
The proof we present here is also similar to the one in \cite{FMT86} (see also \cite{TT88}).
Denote by $\tau (t)= T(t)-(1+\|T_0\|_{\infty} + \|Q\|_{\infty} \; t).$ It is clear that $\tau$
satisfies:
\begin{eqnarray}
&&\hskip-.168in
\frac{\pp \tau}{\pp t}  + v \cdot {\nabla}_{H} \tau + w
\frac{\pp \tau}{\pp z}  +  L_2 \tau =  Q-\|Q\|_{\infty}. \label{TAU}
\end{eqnarray}
Let $\tau^+ = \max\{0, \tau\}$ which belongs to $H^1(\Om)$ and satisfies
\begin{eqnarray}
&&\hskip-.8in  \tau^+(z=0)=\tau^+(z=-h)=0.
\label{TB-1}
\end{eqnarray}
Taking the  inner product of the  equation (\ref{TAU}) with
$\tau^+$ in $L^2(\Om)$ and applying the boundary conditions (\ref{TB-1}), we get
\begin{eqnarray}
&&\hskip-.68in \frac{1}{2} \frac{d \|\tau^+\|_{2}^{2}}{dt}
+ \frac{1}{R_3}\|\pp_z\tau^+\|_2^2
= \int_{\Om} (Q-\|Q\|_{\infty}) \tau^+  \; dxdydz -\int_{\Om} \left[ v
\cdot {\nabla}_{H} \tau +w \,\pp_z \tau
\right]\; \tau^+ \; dxdydz.   \label{TEST-1}
\end{eqnarray}
By integration by parts and using the boundary conditions (\ref{EQ7}) and (\ref{TB-1}), we get
\begin{eqnarray}
&&\hskip-.065in  \int_{\Om} \left[ v
\cdot {\nabla}_{H} \tau +w \,\pp_z \tau
\right]\; \tau^+ \; dxdydz =0.   \label{TEST-2}
\end{eqnarray}
Thus,
\begin{eqnarray}
&&\hskip-.68in \frac{1}{2} \frac{d \|\tau^+\|_{2}^{2}}{dt}
+ \frac{1}{R_3}\|\pp_z\tau^+\|_2^2  \; dxdydz
 = \int_{\Om} (Q-\|Q\|_{\infty}) \tau^+  \; dxdydz \leq 0.   \label{TEST-3}
\end{eqnarray}
Therefore, we obtain
\begin{eqnarray}
&&\hskip-.68in  \|\tau^+(t)\|_{2}^{2} \leq \|\tau^+(t=0)\|_{2}^{2}=0.   \label{TEST-4}
\end{eqnarray}
Thus, $\tau^+(t)\equiv 0.$ As a result, we have
\begin{eqnarray}
&&\hskip-.68in  T(t) \leq 1+\|T_0\|_{\infty} + \|Q\|_{\infty} \; t.   \label{TEST-5}
\end{eqnarray}
By applying similar arguments, we also have
\begin{eqnarray}
&&\hskip-.68in  T(t) \geq -(1+\|T_0\|_{\infty} + \|Q\|_{\infty} \; t).   \label{TEST-6}
\end{eqnarray}
Therefore, $T$ satisfies the following $L^{\infty}-$estimate:
\begin{eqnarray}
&&\hskip-.68in    \|T (t)\|_{\infty}
 \leq K_2 = 1+\|T_0\|_{\infty} + \|Q\|_{\infty} \; t.   \label{K-2}
\end{eqnarray}

\vskip0.1in

\subsection{$\|\widetilde{v} \|_6$ estimates}

Taking the  inner product of the  equation (\ref{EQ4}) with
$|\widetilde{v}|^4 \widetilde{v}$ in $L^2(\Om)$ and using the boundary conditions
(\ref{EQ6}),  we get
\begin{eqnarray*}
&&\hskip-.168in \frac{1}{6} \frac{d \| \widetilde{v} \|_{6}^{6} }{d
t} + \frac{1}{R_1} \int_{\Om} \left(|{\nabla}_{H} \widetilde{v}|^2
|\widetilde{v}|^{4} + \left|{\nabla}_{H} |\widetilde{v}|^2 \right|^2
|\widetilde{v}|^{2} \right) \; dxdydz + \frac{1}{R_2} \int_{\Om}
\left(|\widetilde{v}_z|^2 |\widetilde{v}|^{4}
+ \left|\pp_z |\widetilde{v}|^2 \right|^2 |\widetilde{v}|^{2} \right) \; dxdydz    \\
&&\hskip-.165in = - \int_{\Om} \left\{  (\widetilde{v} \cdot
{\nabla}_{H}) \widetilde{v} - \left( \int_{-h}^z {\nabla}_{H} \cdot
\widetilde{v}(x,y, \xi,t) d\xi \right)  \frac{\pp \widetilde{v}}{\pp
z} +(\widetilde{v} \cdot {\nabla}_{H} ) \overline{v}+ (\overline{v}
\cdot {\nabla}_{H}) \widetilde{v} - \overline{ \left[(\widetilde{v}
\cdot {\nabla}_{H}) \widetilde{v}
+ ({\nabla}_{H} \cdot \widetilde{v}) \; \widetilde{v}\right]} \right.  \\
&&\hskip-.06in
 \left. + f_0 \vec{k} \times \widetilde{v}
 - {\nabla}_{H}  \left( \int_{-h}^z
T(x,y,\xi,t) d\xi -\frac{1}{h} \int_{-h}^0 \int_{-h}^z
T(x,y,\xi,t) d\xi dz  \right)    \right\} \cdot
|\widetilde{v}|^{4} \widetilde{v} \; dxdydz.
\end{eqnarray*}
Observe, again, that
\begin{eqnarray}
&&\hskip-.065in
\left( f_0 \vec{k} \times \widetilde{v} \right)  \cdot
|\widetilde{v}|^{4} \widetilde{v}   =0.
\label{D6-1}
\end{eqnarray}
Moreover, by integration by parts and the boundary conditions
(\ref{EQ6}),  we also get
\begin{eqnarray}
&&\hskip-.065in   - \int_{\Om} \left[ (\widetilde{v} \cdot
{\nabla}_{H}) \widetilde{v} - \left( \int_{-h}^z {\nabla}_{H} \cdot
\widetilde{v}(x,y, \xi,t) d\xi \right)  \frac{\pp \widetilde{v}}{\pp
z} \right] \cdot |\widetilde{v}|^{4} \widetilde{v} \; dxdydz =0.
\label{D6-11}
\end{eqnarray}
Furthermore, by virtue of (\ref{EQ2}) and by the boundary conditions (\ref{EQ6}), in particular the horizontal
periodic boundary conditions, we have
\begin{eqnarray}
&&\hskip-.065in \int_{\Om} (\overline{v} (x, y, t) \cdot
{\nabla}_{H}) \widetilde{v} (x, y, z, t) \cdot |\widetilde{v}(x, y,
z, t)|^{4} \widetilde{v}(x, y, z, t) \;dxdydz =0. \label{D6-3}
\end{eqnarray}
Thus, (\ref{D6-1})--(\ref{D6-3}) imply
\begin{eqnarray}
&&\hskip-.168in \frac{1}{6} \frac{d \| \widetilde{v} \|_{6}^{6} }{d
t} + \frac{1}{R_1} \int_{\Om} \left(|{\nabla}_{H} \widetilde{v}|^2
|\widetilde{v}|^{4} + \left|{\nabla}_{H} |\widetilde{v}|^2 \right|^2
|\widetilde{v}|^{2} \right) \; dxdydz + \frac{1}{R_2} \int_{\Om}
\left(|\widetilde{v}_z|^2 |\widetilde{v}|^{4}
+ \left|\pp_z |\widetilde{v}|^2 \right|^2 |\widetilde{v}|^{2} \right) \; dxdydz
\nonumber  \\
&&\hskip-.165in
 = - \int_{\Om} \left\{
(\widetilde{v} \cdot {\nabla}_{H} ) \overline{v}  - \overline{
(\widetilde{v} \cdot {\nabla}_{H}) \widetilde{v}
+ ({\nabla}_{H} \cdot \widetilde{v}) \; \widetilde{v}}   \right.
\nonumber \\
&&\hskip-.06in
 \left.
- {\nabla}_{H}  \left( \int_{-h}^z T(x,y,\xi,t) d\xi -\frac{1}{h}
\int_{-h}^0 \int_{-h}^z T(x,y,\xi,t) d\xi dz  \right\}    \right)
\cdot |\widetilde{v}|^{4} \widetilde{v} \; dxdydz.
\label{D6-33}
\end{eqnarray}
Notice that by integration by parts and using the boundary conditions
(\ref{EQ6}), in particular the horizontal
periodic boundary conditions, we have
\begin{eqnarray}
&&\hskip-.165in
 - \int_{\Om} \left[
(\widetilde{v} \cdot {\nabla}_{H} ) \overline{v}  - \overline{\left[
(\widetilde{v} \cdot {\nabla}_{H}) \widetilde{v}
+ ({\nabla}_{H} \cdot \widetilde{v}) \; \widetilde{v} \right]}   \right.  \\
&&\hskip-.06in \left. - {\nabla}_{H}  \left( \int_{-h}^z
T(x,y,\xi,t) d\xi -\frac{1}{h} \int_{-h}^0 \int_{-h}^z T(x,y,\xi,t)
d\xi dz  \right)    \right] \cdot |\widetilde{v}|^{4}
\widetilde{v} \; dxdydz
\nonumber \\
&&\hskip-.165in
 = \int_{\Om} \left[
({\nabla}_{H} \cdot \widetilde{v})  \; \overline{v} \cdot
|\widetilde{v}|^{4} \widetilde{v} + (\widetilde{v} \cdot
{\nabla}_{H} ) (|\widetilde{v}|^{4} \widetilde{v}) \cdot
\overline{v}
 - \sum_{k=1}^2 \sum_{j=1}^3 \overline{ \widetilde{v}^k   \widetilde{v}^j} \;
 \pp_{x_k}  (|\widetilde{v}|^{4}
 \widetilde{v}^j) \right. \nonumber \\
&&\hskip-.06in
 \left.
-  \left( \int_{-h}^z T(x,y,\xi,t) d\xi -\frac{1}{h} \int_{-h}^0
\int_{-h}^z T(x,y,\xi,t) d\xi dz  \right) {\nabla}_{H} \cdot
(|\widetilde{v}|^{4} \widetilde{v})   \right]  \; dxdydz.
\label{D6-333}
\end{eqnarray}
As a result, we obtain
\begin{eqnarray}
&&\hskip-.168in \frac{1}{6} \frac{d \| \widetilde{v} \|_{6}^{6} }{d
t} + \frac{1}{R_1} \int_{\Om} \left(|{\nabla}_{H} \widetilde{v}|^2
|\widetilde{v}|^{4} + \left|{\nabla}_{H} |\widetilde{v}|^2 \right|^2
|\widetilde{v}|^{2} \right) \; dxdydz + \frac{1}{R_2} \int_{\Om}
\left(|\widetilde{v}_z|^2 |\widetilde{v}|^{4}
+ \left|\pp_z |\widetilde{v}|^2 \right|^2 |\widetilde{v}|^{2} \right) \; dxdydz
\nonumber \\
&&\hskip-.165in \leq C \int_{M} \left[ |\overline{v}| \int_{-h}^0
|{\nabla}_{H} \widetilde{v}| \; | \widetilde{v}|^5 \; dz \right]
\; dxdy   +C \int_{M} \left[ \left( \int_{-h}^0
|\widetilde{v}|^2 \; dz \right) \left( \int_{-h}^0 |{\nabla}_{H}
\widetilde{v}| \;
| \widetilde{v}|^4 \; dz \right) \right] \; dxdy   \nonumber  \\
&&\hskip-.065in +C \int_{M} \left[ \overline{|T|} \int_{-h}^0
|{\nabla}_{H} \widetilde{v}| \; | \widetilde{v}|^4 \; dz \right]
\; dxdy.
\label{D6-3333}
\end{eqnarray}
Therefore, by the Cauchy--Schwarz inequality and H\"{o}lder inequality
we reach
\begin{eqnarray}
&&\hskip-.168in \frac{1}{6} \frac{d \| \widetilde{v} \|_{6}^{6} }{d
t} + \frac{1}{R_1} \int_{\Om} \left(|{\nabla}_{H} \widetilde{v}|^2
|\widetilde{v}|^{4} + \left|{\nabla}_{H} |\widetilde{v}|^2 \right|^2
|\widetilde{v}|^{2} \right) \; dxdydz + \frac{1}{R_2} \int_{\Om}
\left(|\widetilde{v}_z|^2 |\widetilde{v}|^{4}
+ \left|\pp_z |\widetilde{v}|^2 \right|^2 |\widetilde{v}|^{2} \right) \; dxdydz   \nonumber      \\
&&\hskip-.165in \leq C \int_{M} \left[ |\overline{v}| \left(
\int_{-h}^0 |{\nabla}_{H} \widetilde{v}|^2 \; | \widetilde{v}|^4 \;
dz \right)^{1/2} \left( \int_{-h}^0  | \widetilde{v}|^6 \; dz
\right)^{1/2} \right]
\; dxdy  \nonumber     \\
&&\hskip-.06in +C \int_{M} \left[ \left( \int_{-h}^0
|\widetilde{v}|^2 \; dz \right) \left( \int_{-h}^0 |{\nabla}_{H}
\widetilde{v}|^2 \; | \widetilde{v}|^4 \; dz \right)^{1/2}
\left( \int_{-h}^0  | \widetilde{v}|^4 \; dz \right)^{1/2} \right] \; dxdy  \nonumber     \\
&&\hskip-.06in +C \int_{M} \left[ \overline{|T|} \left( \int_{-h}^0
|{\nabla}_{H} \widetilde{v}|^2 \; | \widetilde{v}|^4 \; dz
\right)^{1/2} \left( \int_{-h}^0  | \widetilde{v}|^4 \; dz
\right)^{1/2} \right]
\; dxdy
\end{eqnarray}
Moreover,
\begin{eqnarray}
&&\hskip-.168in \frac{1}{6} \frac{d \| \widetilde{v} \|_{6}^{6} }{d
t} + \frac{1}{R_1} \int_{\Om} \left(|{\nabla}_{H} \widetilde{v}|^2
|\widetilde{v}|^{4} + \left|{\nabla}_{H} |\widetilde{v}|^2 \right|^2
|\widetilde{v}|^{2} \right) \; dxdydz + \frac{1}{R_2} \int_{\Om}
\left(|\widetilde{v}_z|^2 |\widetilde{v}|^{4}
+ \left|\pp_z |\widetilde{v}|^2 \right|^2 |\widetilde{v}|^{2} \right) \; dxdydz   \nonumber      \\
&&\hskip-.165in \leq C \|\overline{v}\|_{L^4(M)}
 \left( \int_{\Om}  |{\nabla}_{H} \widetilde{v}|^2 \; | \widetilde{v}|^4 \; dxdydz \right)^{1/2}
\left( \int_M \left( \int_{-h}^0  | \widetilde{v}|^6 \; dz
\right)^{2}
\; dxdy \right)^{1/4}  \nonumber      \\
&&\hskip-.06in +C \left( \int_M \left( \int_{-h}^0  |
\widetilde{v}|^2 \; dz \right)^{4} \; dxdy \right)^{1/4} \left(
\int_{\Om}  |{\nabla}_{H} \widetilde{v}|^2 \; | \widetilde{v}|^4 \;
dxdydz \right)^{1/2} \left( \int_M \left( \int_{-h}^0  |
\widetilde{v}|^4 \; dz \right)^{2}
\; dxdy \right)^{1/4}  \nonumber      \\
&&\hskip-.06in +C \|\,\overline{|T|}\,\|_{L^4(M)} \left( \int_{\Om}
|{\nabla}_{H} \widetilde{v}|^2 \; | \widetilde{v}|^4 \; dxdydz
\right)^{1/2} \left( \int_M \left( \int_{-h}^0  | \widetilde{v}|^4
\; dz \right)^{2} \; dxdy \right)^{1/4}. \label{D6_22}
\end{eqnarray}
Using the Minkowsky  inequality (\ref{MKY})  we get
\begin{eqnarray}
&&\hskip-.168in \left( \int_M \left( \int_{-h}^0  |
\widetilde{v}|^6 \; dz \right)^{2} \; dxdy \right)^{1/2} \leq  C
\int_{-h}^0 \left( \int_M  | \widetilde{v}|^{12} \; dxdy
\right)^{1/2}  \; dz.   \label{D6_2}
\end{eqnarray}
Thanks to (\ref{TWE}),
\begin{eqnarray}
&&\hskip-.168in \int_M  | \widetilde{v}|^{12}  \; dxdy \leq C_0
\left( \int_M | \widetilde{v}|^{6}  \; dxdy \right) \left( \int_M |
\widetilde{v}|^{4} |{\nabla}_{H} \widetilde{v} |^2  \; dxdy \right)
+ \left( \int_M | \widetilde{v}|^{6}  \; dxdy \right)^2.
\label{D6_3}
\end{eqnarray}
Thus, by the Cauchy--Schwarz inequality we obtain
\begin{eqnarray}
&&\hskip-.168in \left( \int_M \left( \int_{-h}^0  | \widetilde{v}|^6
\; dz \right)^{2} \; dxdy \right)^{1/2} \leq  C \|
\widetilde{v}\|_{L^6(\Om)}^{3}  \left( \int_{\Om}
|\widetilde{v}|^{4} |{\nabla}_{H} \widetilde{v} |^2 \; dxdydz
\right)^{1/2} + \| \widetilde{v}\|_{L^6(\Om)}^{6}. \label{M1}
\end{eqnarray}
Similarly, by (\ref{MKY}) and (\ref{SI-2}), we also obtain
\begin{eqnarray}
&&\hskip-.168in \left( \int_M \left( \int_{-h}^0  |
\widetilde{v}|^4 \; dz \right)^{2} \; dxdy \right)^{1/2} \leq  C
\int_{-h}^0 \left( \int_M  | \widetilde{v}|^{8}
\; dxdy \right)^{1/2}  \; dz   \nonumber  \\
&&\hskip-.168in \leq  C  \int_{-h}^0  \|
\widetilde{v}\|_{L^6(M)}^{3} \left( \|{\nabla}_{H}
\widetilde{v}\|_{L^2(M)}  + \| \widetilde{v}\|_{L^2(M)} \right) \;
dz \leq C  \| \widetilde{v}\|_6^{3} \left( \|{\nabla}_{H}
\widetilde{v}\|_2 + \| \widetilde{v}\|_2 \right),   \label{M2}
\end{eqnarray}
and
\begin{eqnarray}
&&\hskip-.168in \left( \int_M \left( \int_{-h}^0  |
\widetilde{v}|^2 \; dz \right)^{4} \; dxdy \right)^{1/4} \leq  C
\int_{-h}^0 \left( \int_M  | \widetilde{v}|^{8}
\; dxdy \right)^{1/4}  \; dz    \nonumber  \\
&&\hskip-.168in \leq  C  \int_{-h}^0  \|
\widetilde{v}\|_{L^6(M)}^{3/2} \left(
 \|{\nabla}_{H} \widetilde{v}\|_{L^2(M)}^{1/2}   + \| \widetilde{v}\|_{L^2(M)}^{1/2} \right) \; dz
\leq C  \| \widetilde{v}\|_6^{3/2} \left( \|{\nabla}_{H}
\widetilde{v}\|_2  + \| \widetilde{v}\|_2 \right)^{1/2}.
\label{M3}
\end{eqnarray}
Therefore, using (\ref{M1})--(\ref{M3}) and (\ref{SI-1}), we reach to
\begin{eqnarray*}
&&\hskip-.168in \frac{1}{6} \frac{d \| \widetilde{v} \|_{6}^{6} }{d
t} + \frac{1}{R_1} \int_{\Om} \left(|{\nabla}_{H} \widetilde{v}|^2
|\widetilde{v}|^{4} + \left|{\nabla}_{H} |\widetilde{v}|^2 \right|^2
|\widetilde{v}|^{2} \right) \; dxdydz + \frac{1}{R_2} \int_{\Om}
\left(|\widetilde{v}_z|^2 |\widetilde{v}|^{4}
+ \left|\pp_z |\widetilde{v}|^2 \right|^2 |\widetilde{v}|^{2} \right) \; dxdydz    \\
&&\hskip-.165in \leq C \|\overline{v}\|_2^{1/2} \; \|{\nabla}_{H}
\overline{v}\|_2^{1/2} \| \widetilde{v}\|_6^{3/2}
 \left( \int_{\Om}  |{\nabla}_{H} \widetilde{v}|^2 \; | \widetilde{v}|^4 \; dxdydz \right)^{3/4}
 +C \| \widetilde{v}\|_6^{3} \left( \|{\nabla}_{H}
\widetilde{v}\|_2 + \| \widetilde{v}\|_2 \right) \left( \int_{\Om}
|{\nabla}_{H} \widetilde{v}|^2 \; | \widetilde{v}|^4 \; dxdydz
\right)^{1/2}
   \\
&&\hskip-.065in +C \|\overline{v}\|_2^{1/2} \; \|{\nabla}_{H} \overline{v}\|_2^{1/2}
\| \widetilde{v}\|_6^{6} +C \|T\|_{\infty} \;  \| \widetilde{v}\|_6^{3/2}
\left( \|{\nabla}_{H} \widetilde{v}\|_2^{1/2} + \|
\widetilde{v}\|_2^{1/2} \right) \left( \int_{\Om}  |{\nabla}_{H}
\widetilde{v}|^2 \; | \widetilde{v}|^4 \; dxdydz \right)^{1/2}.
\end{eqnarray*}
By Young's inequality and Cauchy--Schwarz inequality we have
\begin{eqnarray*}
&&\hskip-.168in
 \frac{d \| \widetilde{v} \|_{6}^{6} }{d t} +
\frac{1}{R_1} \int_{\Om} \left(|{\nabla}_{H} \widetilde{v}|^2
|\widetilde{v}|^{4} + \left|{\nabla}_{H} |\widetilde{v}|^2 \right|^2
|\widetilde{v}|^{2} \right) \; dxdydz + \frac{1}{R_2} \int_{\Om}
\left(|\widetilde{v}_z|^2 |\widetilde{v}|^{4}
+ \left|\pp_z |\widetilde{v}|^2 \right|^2 |\widetilde{v}|^{2} \right) \; dxdydz    \\
&&\hskip-.165in \leq C \|\overline{v}\|_2^{2} \; \|{\nabla}_{H}
\overline{v}\|_2^{2}   \| \widetilde{v}\|_6^{6} +C \|
\widetilde{v}\|_6^{6} \|{\nabla}_{H} \widetilde{v}\|_2^2 +C
\|T\|_{\infty}^4  +C \| \widetilde{v}\|_2^2\| \widetilde{v}\|_6^{6}.
\end{eqnarray*}
Thanks to (\ref{K-1}), (\ref{K-2}) and Gronwall inequality, we get
\begin{eqnarray}
&&\hskip-.68in \| \widetilde{v} (t)\|^6_6 + \int_0^t \left(
\frac{1}{R_1} \int_{\Om} |{\nabla}_{H} \widetilde{v}|^2
|\widetilde{v}|^{4}
 \; dxdydz +
\frac{1}{R_2} \int_{\Om} |\widetilde{v}_z|^2
|\widetilde{v}_z|^{4}  \; dxdydz \right) \leq K_3, \label{K-3}
\end{eqnarray}
where
\begin{eqnarray}
&&\hskip-.68in K_3 = e^{K_1^2 t} \left[ \|v_0\|_{H^1(\Om)}^6 +
K_2^4 \; t \right].   \label{K3}
\end{eqnarray}

\subsection{$\|{\nabla}_{H} \overline{v}\|_2$ estimates}

By taking the inner product of equation (\ref{EQ1}) with $-
{\Dd}_{H} \overline{v}$ in $L^2(M)$, and applying (\ref{EQ2}),
and using the boundary conditions (\ref{EQ6}), we
reach
\begin{eqnarray}
&&\hskip-.68in \frac{1}{2} \frac{d \| {\nabla}_{H} \overline{v}
\|_2^2 }{d t} + \frac{1}{R_1} \|{\Dd}_{H} \overline{v}\|_2^2  =
\int_{M} \left\{ (\overline{v} \cdot {\nabla}_{H} ) \overline{v} +
\overline{ \left[ (\widetilde{v} \cdot {\nabla}_{H}) \widetilde{v} +
({\nabla}_{H} \cdot \widetilde{v}) \; \widetilde{v}\right]} +f_0
\vec{k}\times  \overline{v}\right\} \cdot {\Dd}_{H} \overline{v} \;
dxdy.    \label{DH-1}
\end{eqnarray}
Following the situation for the $2D$ Navier--Stokes equations (cf. e.g.,
\cite{CF88}, \cite{SALMON}) we have
\begin{eqnarray}
&&\hskip-.68in \left| \int_{M}  (\overline{v} \cdot {\nabla}_{H} )
\overline{v} \cdot {\Dd}_{H} \overline{v} \; dxdy \right|
 \leq C \|\overline{v}\|_2^{1/2} \|{\nabla}_{H} \overline{v}\|_2
\; \|{\Dd}_{H}  \overline{v}\|_2^{3/2}.    \label{DH-2}
\end{eqnarray}
By the Cauchy--Schwarz and the H\"{o}lder inequalities, we have
\begin{eqnarray}
&&\hskip-.68in \left|\int_{M} \; \overline{ (\widetilde{v} \cdot
{\nabla}_{H}) \widetilde{v}  + ({\nabla}_{H} \cdot \widetilde{v}) \;
\widetilde{v}}  \cdot {\Dd}_{H} \overline{v} \; dxdy \right| \leq C
\int_M \int_{-h}^0 |\widetilde{v}| \; |{\nabla}_{H} \widetilde{v}|
\; dz \; |{\Dd}_{H} \overline{v}| \; dxdy    \nonumber
\\
&&\hskip-.68in \leq C \int_M \left[ \left( \int_{-h}^0
|\widetilde{v}|^2 \; |{\nabla}_{H} \widetilde{v}| \; dz
\right)^{1/2} \left( \int_{-h}^0  |{\nabla}_{H} \widetilde{v}| \; dz
\right)^{1/2} \; |{\Dd}_{H} \overline{v}| \right] \; dxdy   \nonumber
\\
&&\hskip-.68in \leq C \left[ \int_M  \left( \int_{-h}^0
|\widetilde{v}|^2 \; |{\nabla}_{H} \widetilde{v}| \; dz \right)^{2}
\; dxdy \right]^{1/4} \; \left[ \int_M \left( \int_{-h}^0
|{\nabla}_{H} \widetilde{v}| \; dz \right)^{2} \; dxdy \right]^{1/4}
\; \left[ \int_M
 |{\Dd}_{H} \overline{v}|^2  \; dxdy  \right]^{1/2} \nonumber  \\
&&\hskip-.68in \leq C \|{\nabla}_{H} \widetilde{v}\|_2^{1/2} \left(
\int_{\Om} |\widetilde{v}|^4 |{\nabla}_{H} \widetilde{v}|^2 \;
dxdydz \right)^{1/4}\|{\Dd}_{H} \overline{v}\|_2,   \label{DH-3}
\end{eqnarray}
and
\begin{eqnarray}
&&\hskip-.68in \left|\int_{M} \; f_0 \times \overline{v} \cdot
{\Dd}_{H} \overline{v} \; dxdy \right| \leq C \| \overline{v}\|_2
\|{\Dd}_{H} \overline{v}\|_2.   \label{DH-4}
\end{eqnarray}
Thus, by Young's inequality and the Cauchy--Schwarz inequality, we
have
\begin{eqnarray}
&&\hskip-.68in \frac{d \| {\nabla}_{H} \overline{v} \|_2^2 }{d t} +
\frac{1}{R_1} \|{\Dd}_{H} \overline{v}\|_2^2 \leq C
\|\overline{v}\|_2^2 \|{\nabla}_{H} \overline{v}\|_2^4 + C
\|{\nabla}_{H} \widetilde{v}\|_2^{2} + C \int_{\Om}
|\widetilde{v}|^4 |{\nabla}_{H} \widetilde{v}|^2 \; dxdydz + C
\|\overline{v}\|_2^2.   \label{DH-5}
\end{eqnarray}
By (\ref{K-1}), (\ref{K-3}) and thanks to Gronwall inequality we
obtain
\begin{eqnarray}
&&\hskip-.68in \| {\nabla}_{H} \overline{v} \|_2^2 + \frac{1}{R_1}
\int_0^t |{\Dd}_{H} \overline{v}|_2^2 \; ds \leq K_4,   \label{K-4}
\end{eqnarray}
where
\begin{eqnarray}
&&\hskip-.68in K_4 = e^{K_2^2 t} \left[ \|v_0\|_{H^1(\Om)}^2 + K_2
+K_3 \right].   \label{K4}
\end{eqnarray}

\subsection{$\|v_z\|_6$ estimates}
The {\it a priori}  estimates (\ref{K-1})--(\ref{K-4}) are essentially similar to those obtained in \cite{CT05}.
From now on, we will get new
{\it a priori} estimates of various norms.

Denote by $u=v_z.$ It is clear that $u$ satisfies
\begin{eqnarray}
&&\hskip-.68in
 \frac{\pp u }{\pp t} + L_1 u +
(v \cdot {\nabla}_{H}) u - \left( \int_{-h}^z {\nabla}_{H} \cdot v
(x,y, \xi,t)
d\xi \right) \frac{\pp u}{\pp z}    \nonumber    \\
&&\hskip-.46in + (u \cdot {\nabla}_{H} ) v - ({\nabla}_{H} \cdot v) u
 + f_0 \vec{k} \times u - {\nabla}_{H} T
  = 0. \label{UU}  \\
&&\hskip-.68in \left.   u \right|_{z=0}=  \left. u \right|_{z=-h}=
0.  \label{UU-B}
\end{eqnarray}
Taking the  inner product of the  equation
(\ref{UU}) with $u |u|^4$ in $L^2$,  we get
\begin{eqnarray}
&&\hskip-.68in \frac{1}{6} \frac{d \| u \|_6^6 }{d t} +
\frac{5}{R_1} \| \; |u|^2 \; |{\nabla}_{H} u| \; \|_2^2 +
\frac{5}{R_2}  \|\; |u|^2 \; |\pp_z u|\; \|_2^2     \nonumber  \\
&&\hskip-.65in =- \int_{\Om} \left(  (v \cdot {\nabla}_{H}) u -
\left( \int_{-h}^z {\nabla}_{H} \cdot v (x,y, \xi,t)
d\xi \right) \frac{\pp u}{\pp z}  \right) \cdot u |u|^4 \; dxdydz  \nonumber  \\
&&\hskip-.58in - \int_{\Om} \left( (u \cdot {\nabla}_{H} ) v -
({\nabla}_{H} \cdot v) u + f_0 \vec{k} \times u -{\nabla}_{H} T
\right) \cdot u |u|^4 \; dxdydz.   \label{DZ_1}
\end{eqnarray}
Notice again that
\begin{eqnarray}
&&\hskip-.065in    f_0 \vec{k} \times u  \cdot u |u|^4 =0.
\label{DZ-1}
\end{eqnarray}
 Integrating by parts and using the boundary conditions,
 in particular (\ref{UU-B}),  give
\begin{eqnarray}
&&\hskip-.065in   - \int_{\Om} \left( (v \cdot {\nabla}_{H}) u -
\left( \int_{-h}^z {\nabla}_{H} \cdot v(x,y, \xi,t) d\xi \right)
\frac{\pp u}{\pp  z} \right) \cdot u |u|^4\; dxdydz =0.
\label{DZ-11}
\end{eqnarray}
Thus, by (\ref{DZ-1}), (\ref{DZ-11}) and H\"{o}lder inequality,  we
have
\begin{eqnarray}
&&\hskip-.68in \frac{1}{6} \frac{d \| u \|_6^6 }{d t} +
\frac{5}{R_1} \| \; |u|^2 \; |{\nabla}_{H} u| \; \|_2^2 +
\frac{5}{R_2}  \|\; |u|^2 \; |\pp_z u|\; \|_2^2     \nonumber    \\
&&\hskip-.65in =- \int_{\Om} \left(  (u \cdot {\nabla}_{H} ) v -
({\nabla}_{H} \cdot v) u
-{\nabla}_{H} T  \right) \cdot u |u|^4\; dxdydz    \nonumber   \\
&&\hskip-.65in \leq C \int_{\Om}  |v|  \,|u|^5 \, |{\nabla}_{H} u|
\; dxdydz +C \int_{\Om} |T|   \,|u|^4 \,
|{\nabla}_{H} u|  \; dxdydz   \nonumber   \\
&&\hskip-.65in \leq C  \left( \|v\|_6 \|u^3\|_3 \|{\nabla}_{H}
u^2\|_2 +
\|T\|_6 \; \| u^2 \|_3  \| {\nabla}_{H} u^3  \|_2 \right)  \nonumber  \\
&&\hskip-.65in \leq C  \left( \|v\|_6 \|u\|_6^{3/2} \|{\nabla}_{H}
u^3\|_2^{3/2} + \|T\|_6 \; \| u \|_6^2  \| {\nabla}_{H} u^3  \|_2
\right). \label{DZ_2}
\end{eqnarray}
Thanks to the Cauchy--Schwarz inequality, we have
\begin{eqnarray}
&&\hskip-.68in \frac{d \| u \|_6^6 }{d t}  + \frac{1}{R_1} \| \;
|u|^2 \; |{\nabla}_{H} u| \; \|_2^2 +
\frac{1}{R_2}  \|\; |u|^2 \; |\pp_z u|\; \|_2^2         \nonumber   \\
&&\hskip-.65in  \leq C \left( 1+\|v\|_6^4  \right) \; \| u \|_6^6 + \|T\|_6^6  \\
&&\hskip-.65in  \leq C \left(1+ \|{\nabla}_{H} \overline{v}\|_2^4
+\|\widetilde{v}\|_6^4  \right) \; \| u \|_6^6  + \|T\|_6^6.   \label{DZ_3}
\end{eqnarray}
Using  (\ref{K-1}), (\ref{K-3}), (\ref{K-4}), and Gronwall
inequality, we get
\begin{eqnarray}
&&\hskip-.68in \| u \|_6^6 +  \int_0^t \left[ \; \frac{1}{R_1} \| \;
|u|^2 \; |{\nabla}_{H} u| \; \|_2^2 + \frac{1}{R_2}  \|\; |u|^2 \;
|\pp_z u|\; \|_2^2  \right]  \; ds \leq K_5,   \label{K-5}
\end{eqnarray}
where
\begin{eqnarray}
&&\hskip-.68in K_5 = e^{(1+K_3^{2/3}+K_4^{2}) t} \left[ \|\pp_z
v_0\|_{H^1(\Om)}^6 + K_2^6 \; t \right].   \label{K5}
\end{eqnarray}

\subsection{$\|v_{zz}\|_2$ estimates}

Taking the  inner product of the  equation (\ref{UU}) with $-
u_{zz}$ in $L^2(\Om)$ and recalling that $u=v_z,$ which satisfies the boundary
condition (\ref{UU-B}), we get
\begin{eqnarray*}
&&\hskip-.168in \frac{1}{2} \frac{d \| u_z \|_2^2 }{d t} +
\frac{1}{R_1} \|  {\nabla}_{H} u_z \|_2^2 +
\frac{1}{R_2}  \|u_{zz} \|_2^2     \\
&&\hskip-.065in =\int_{\Om} \left(  (v \cdot {\nabla}_{H}) u -
\left( \int_{-h}^z {\nabla}_{H} \cdot v (x,y, \xi,t)
d\xi \right) \frac{\pp u}{\pp z}  \right) \cdot u_{zz}  \; dxdydz \\
&&\hskip-.058in + \int_{\Om} \left( (u \cdot {\nabla}_{H} ) v -
({\nabla}_{H} \cdot v) u  + f_0 \vec{k} \times u -{\nabla}_{H} T
\right) \cdot u_{zz}\;
dxdydz    \\
&&\hskip-.065in = -\int_{\Om}   \left[ (u \cdot {\nabla}_{H}) u +
 (u_z \cdot {\nabla}_{H} ) v +(u \cdot {\nabla}_{H} ) u - ({\nabla}_{H}
\cdot u) u - ({\nabla}_{H} \cdot v) u_z  \right] \cdot u_{z}\;
dxdydz - \int_{\Om}  T_z ({\nabla}_{H} \cdot u_z) \;
dxdydz  \\
&&\hskip-.065in \leq C \|u\|_6 \|{\nabla}_{H} u \|_2 \|u_z\|_3 +C
\|v\|_6 \|{\nabla}_{H} u_z \|_2 \|u_z\|_3
 + \|T_z\|_2  \| {\nabla}_{H} u_z \|_2    \\
&&\hskip-.065in \leq C \left[ \|u\|_6 \|{\nabla}_{H} u \|_2 +
\|v\|_6 \|{\nabla}_{H} u_z \|_2  \right] \|u_z\|_2^{1/2} \left(
\|{\nabla}_{H} u_z\|_2^{1/2} +\| u_{zz} \|_2^{1/2}\right)
 + \|T_z\|_2  \| {\nabla}_{H} u_z \|_2.
\end{eqnarray*}
By the Cauchy--Schwarz  and Young's inequalities, we have
\begin{eqnarray*}
&&\hskip-.68in
 \frac{1}{2} \frac{d \| u_z \|_2^2 }{d t} +
\frac{1}{R_1} \|  {\nabla}_{H} u_z \|_2^2 +
\frac{1}{R_2}  \|u_{zz} \|_2^2          \\
&&\hskip-.65in  \leq C \left( \|v\|_6^4 + \|u\|_6^4 \right) \; \|
u_z\|_2^2 + C \|{\nabla}_{H} u \|_2^2 + C \|T_z\|_2^2.
\end{eqnarray*}
Applying  (\ref{K-1}), (\ref{K-3}), (\ref{K-5}),
and Gronwall inequality yield
\begin{eqnarray}
&&\hskip-.68in \| u_z \|_2^2+  \int_0^t \left[ \; \frac{1}{R_1} \|
{\nabla}_{H} u_z \|_2^2 + \frac{1}{R_2}  \|u_{zz} \|_2^2   \right]
\; ds \leq K_6,   \label{K-6}
\end{eqnarray}
where
\begin{eqnarray}
&&\hskip-.68in K_6 = C e^{(K_3^{2/3}+K_5^{2/3}) t} \left[
\|v_0\|_{H^1(\Om)}^2 + K_1 \right].   \label{K6}
\end{eqnarray}

\subsection{$\|{\nabla}_{H} \times v_z\|_2^2 + \| {\nabla}_{H} \cdot v_z + R_1 T
\|_2^2$ estimates}

Let $\beta$ be the solution of the following two--dimensional
elliptic problem with periodic boundary conditions:
\begin{eqnarray}
{\Dd}_{H} \beta = {\nabla}_{H} T,  \quad  \int_M \beta \; dxdy =0,   \label{BETA-1}
\end{eqnarray}
where $z$ is considered as a parameter. Roughly speaking, $\beta$ is like the potential
vorticity. Notice that
\begin{eqnarray}
{\nabla}_{H} \cdot \beta =T, \qquad {\nabla}_{H} \times \beta =0.   \label{BETA-2}
\end{eqnarray}
Recall that $u=v_z.$ We denote by
\begin{eqnarray}
&& \zeta=u+R_1 \beta,  \label{ZZZ}  \\
&& \eta=\left({\nabla}_H^{\perp
} \cdot \zeta \right)= {\nabla}_H^{\perp} \cdot u
= \pp_x u_2 - \pp_y u_1,  \label{EEE}  \\
&& \tt = \left( {\nabla}_{H} \cdot \zeta \right) = {\nabla}_{H}
\cdot u + R_1 T = \pp_x u_1 + \pp_y u_2 + R_1 T. \label{TTT}
\end{eqnarray}
Applying (\ref{DIV-CUR}) for $\zeta$ and
$\beta$ with $m \geq 0, 1< q <\infty$,  using (\ref{BETA-2})--(\ref{TTT}), we obtain
\begin{eqnarray}
&& \|{\nabla}_{H} u\|_{W^{m,q}(M)} \leq C\left( \|{\nabla}_{H}
\zeta\|_{W^{m,q}(M)}
+ R_1 \|{\nabla}_{H} \beta\|_{W^{m,q}(M)}  \right) \nonumber   \\
&& \leq C \left( \|\eta \|_{W^{m,q}(M)} +\|\tt \|_{W^{m,q}(M)}
+ \| T \|_{W^{m,q}(M)}  \right), \label{INQ}
\end{eqnarray}
where, again, we consider $z$ as a parameter;
consequently, the constant $C$ is independent of
$z$. By applying the operator ${\nabla}_H^{\perp}
\cdot$ to equation (\ref{UU}) we obtain
\begin{eqnarray}
&&\hskip-.68in
 \frac{\pp \eta }{\pp t} + L_1 \eta  + {\nabla}_H^{\perp} \cdot \left[
(v \cdot {\nabla}_{H}) u - \left( \int_{-h}^z {\nabla}_{H} \cdot v
(x,y, \xi,t) d\xi \right) \frac{\pp u}{\pp z}  +  (u \cdot
{\nabla}_{H} ) v - ({\nabla}_{H} \cdot v) u \right] \nonumber  \\
&&\hskip-.6in  - f_0 \left(R_1\, T- \tt\right)
 =0. \label{ETA}
\end{eqnarray}
Then,  multiplying
 equation (\ref{EQ5}) by $R_1$ and adding to the above equation we reach
\begin{eqnarray}
&&\hskip-.68in
 \frac{\pp \tt }{\pp t} + L_1 \tt + {\nabla}_{H} \cdot \left[
(v \cdot {\nabla}_{H}) u - \left( \int_{-h}^z {\nabla}_{H} \cdot v
(x,y, \xi,t) d\xi \right) \frac{\pp u}{\pp z}  +  (u \cdot
{\nabla}_{H} ) v -
({\nabla}_{H} \cdot v) u  \right] - f_0\, \eta \nonumber   \\
&&\hskip-.6in   + R_1\left[ v \cdot {\nabla}_{H} T - \left(
\int_{-h}^z {\nabla}_{H}
 \cdot v (x,y, \xi,t) d\xi \right) \left( \frac{\pp T}{\pp
 z}+\frac{1}{h}\right)\right]
 =R_1Q +R_1\, \left(\frac{1}{R_3}-\frac{1}{R_2}\right) T_{zz}. \label{TT}
\end{eqnarray}
Taking the  inner product of  equation (\ref{ETA}) with $\eta$
in $L^2(\Om)$ and the  equation (\ref{TT}) with $\tt$ in $L^2(\Om)$, integrating by parts and
observing that $\left.   \eta \right|_{z=0}=  \left. \eta
\right|_{z=-h}=\left.   \tt \right|_{z=0}=  \left. \tt
\right|_{z=-h}= 0$, we get
\begin{eqnarray*}
&&\hskip-.168in \frac{1}{2} \frac{d \left( \| \eta \|_2^2+ \| \tt
\|_2^2 \right)  }{d t} + \frac{1}{R_1} \left( \|  {\nabla}_{H} \eta
\|_2^2+ \|  {\nabla}_{H} \tt  \|_2^2 \right)   +
\frac{1}{R_2} \left(  \|\pp_z \eta \|_2^2  + \|\pp_z \tt \|_2^2 \right)  \\
&&\hskip-.158in  = - \int_{\Om} \left\{ {\nabla}_H^{\perp} \cdot \left[
(v \cdot {\nabla}_{H}) u - \left( \int_{-h}^z {\nabla}_{H} \cdot v
(x,y, \xi,t) d\xi \right) \frac{\pp u}{\pp z}  +  (u \cdot
{\nabla}_{H} ) v - ({\nabla}_{H} \cdot v) u \right] \; \eta -f_0 \,
R_1 T \, \eta\right\}
\; dxdydz   \\
&&\hskip-.05in - \int_{\Om} \left\{ {\nabla}_{H} \cdot \left[ (v
\cdot {\nabla}_{H}) u - \left( \int_{-h}^z {\nabla}_{H} \cdot v
(x,y, \xi,t) d\xi \right) \frac{\pp u}{\pp z}  +  (u \cdot
{\nabla}_{H} ) v - ({\nabla}_{H} \cdot v) u \right]\, \tt\right\}
\; dxdydz   \\
&&\hskip-.05in - \int_{\Om} \left\{ R_1\left[ v \cdot {\nabla}_{H} T
- \left( \int_{-h}^z {\nabla}_{H}
 \cdot v (x,y, \xi,t) d\xi \right) \left( \frac{\pp T}{\pp
 z}+\frac{1}{h}\right)
 - Q +\, \left(\frac{1}{R_3}-\frac{1}{R_2}\right) T_{zz} \right] \, \tt\right\}
\; dxdydz
   \\
&&\hskip-.165in \leq C \|Q\|_2 \|\tt\|_2 +  C \int_{\Om} \left(|v|
\,|{\nabla}_{H} u| \,+|u| \,|{\nabla}_{H} v| +|u_z|
\,\int_{-h}^0|{\nabla}_{H} \cdot v| \; d\xi \right)
\left(|{\nabla}_{H} \eta|+|{\nabla}_{H}
\tt|\right)   \; dxdydz +C \|T\|_2 \|\eta\|_2   \\
&&\hskip-.06in  + C \int_{\Om}  \left[ |{\nabla}_{H} v| |T| |\tt|
+|v| |T| |{\nabla}_{H} \tt| +  |T_z| \; \left( \int_{-h}^0
|{\nabla}_{H} \cdot v| \; d\xi\right) \;|\tt|  \right] \; dxdydz +C
\|v\|_2 \|{\nabla}_{H} \tt\|_2  + C \|T_z\|_2 \|\tt_z\|_2
\\
&&\hskip-.165in \leq C \|Q\|_2   \|\tt\|_2 +C \int_{\Om} \left(|v|
\,|{\nabla}_{H} u| \,+|u| \,|{\nabla}_{H} v| +|u_z|
\,\int_{-h}^0(|\tt|+ |T|) \; d\xi \right) \left(|{\nabla}_{H}
\eta|+|{\nabla}_{H} \tt|\right)
\; dxdydz   +C \|T\|_2 \|\eta\|_2 \\
&&\hskip-.06in  + C \int_{\Om}  \left[ |{\nabla}_{H} v| |T| |\tt| +
|v| |T| |{\nabla}_{H} \tt| +  |T_z| \; \left( \int_{-h}^0
(|\tt|+|T|) \; d\xi \right) | \tt| \right] \; dxdydz +C \|v\|_2
\|{\nabla}_{H} \tt\|_2 + C \|T_z\|_2 \|\tt_z\|_2.
\end{eqnarray*}
Using H\"{o}lder inequality, and inequalities (\ref{MAIN-1}) and
(\ref{INQ}), we obtain
\begin{eqnarray*}
&&\hskip-.168in \frac{1}{2} \frac{d \left( \| \eta \|_2^2+ \| \tt
\|_2^2 \right)  }{d t} + \frac{1}{R_1} \left( \|  {\nabla}_{H} \eta
\|_2^2+ \|  {\nabla}_{H} \tt  \|_2^2 \right)   + \frac{1}{R_2}
\left( \|\pp_z \eta \|_2^2  + \|\pp_z \tt \|_2^2 \right)
 \\
&&\hskip-.165in \leq C \left(\|v\|_6  \| {\nabla}_{H} u \|_3 +
\|u\|_6 \| {\nabla}_{H} v \|_3  + \|u_z\|_2^{\frac{1}{2}}  \|
{\nabla}_{H} u_z \|_2^{\frac{1}{2}} \| \tt \|_2^{\frac{1}{2}}  \|
{\nabla}_{H} \tt \|_2^{\frac{1}{2}} + \|T\|_{\infty}  \| u_z \|_2
\right) \left( \| {\nabla}_{H}
\eta \|_2 +\|  {\nabla}_{H} \tt \|_2\right)  \\
&&\hskip-.06in +C \|T\|_2 \|\eta\|_2 +C \|{\nabla}_{H} v\|_2 \|T
\|_{\infty} \|\tt\|_2 + C \|v\|_2 \|T \|_{\infty} \|{\nabla}_{H}
\tt\|_2
  + C \|T_z\|_2  \| \tt \|_2 \| {\nabla}_{H} \tt \|_2
+ C \|T\|_{\infty}  \| T_z \|_2  \| \tt  \|_2 \\
&&\hskip-.06in
+C \|Q\|_2 \|\tt\|_2
+C \|v\|_2 \|{\nabla}_{H} \tt\|_2  + C \|T_z\|_2 \|\tt_z\|_2 \\
&&\hskip-.165in \leq C  (\|T\|_3 + \| \eta \|_2^{\frac{1}{2}} \|
{\nabla}_{H} \eta \|_2^{\frac{1}{2}}+ \| \tt \|_2^{\frac{1}{2}} \|
{\nabla}_{H} \tt \|_2^{\frac{1}{2}})  \; \left(\|v\|_6 +\|u\|_6
\right) \;
  \left(\| {\nabla}_{H} \eta \|_2 +\| {\nabla}_{H} \tt\|_2\right) +    \\
&&\hskip-.06in  +C \left( \|u_z\|_2^{1/2}  \| {\nabla}_{H} u_z
\|_2^{1/2} \| \tt \|_2^{1/2}  \|  {\nabla}_{H} \tt  \|_2^{1/2}    +
\|T\|_{\infty} \| u_z \|_2 + \|T_z\|_2  \| \tt \|_2 +\|v\|_2 \|T
\|_{\infty} \right) \left(\| {\nabla}_{H} \eta \|_2
+\| {\nabla}_{H} \tt\|_2\right)   \\
&&\hskip-.06in  +C \|T\|_2 \|\eta\|_2 + C \left( \|{\nabla}_{H} v
\|_2 \; \|T \|_{\infty} + \|T\|_{\infty}  \| T_z \|_2  \right) \|
\tt \|_2+C \|Q\|_2 \|\tt\|_2+C \|v\|_2 \|{\nabla}_{H} \tt\|_2  + C
\|T_z\|_2 \|\tt_z\|_2.
\end{eqnarray*}
By Young's and the Cauchy--Schwarz inequalities, we have
\begin{eqnarray*}
&&\hskip-.68in  \frac{d \left( \| \eta \|_2^2+ \| \tt \|_2^2 \right)
}{d t} + \frac{1}{R_1} \left( \|  {\nabla}_{H} \eta \|_2^2+ \|
{\nabla}_{H} \tt  \|_2^2 \right)   +
\frac{1}{R_2} \left(  \|\pp_z \eta \|_2^2  + \|\pp_z \tt \|_2^2 \right)       \\
&&\hskip-.65in  \leq C \left(1+ \|v\|_6^4 + \|u\|_6^4+ \|T_z\|_2^2 +
\|u_z\|_2^2\|{\nabla}_{H} u_z\|_2^2 \right) \; \left(
\| \eta \|_2^2+ \| \tt \|_2^2 \right) \\
&&\hskip-.5in
+ C \|T\|_2^2 + C \|Q\|_2^2 +C \left( \|v\|_6^2 +
\|u\|_6^2 + \|u_z\|_2^2 +\|v\|_2^2+\|{\nabla}_{H}
v\|_2^2+\|T_z\|_2^2 \right) \left( 1+\|T\|_{\infty}^2\right).
\end{eqnarray*}
Thanks to (\ref{K-1}), (\ref{K-2}), (\ref{K-3}), (\ref{K-5}),
(\ref{K-6}), and Gronwall inequality, we have
\begin{eqnarray}
&&\hskip-.68in \| \eta \|_2^2+ \| \tt \|_2^2 +
 \int_0^t \left[ \frac{1}{R_1} \left( \|  {\nabla}_{H} \eta  \|_2^2+ \|  {\nabla}_{H} \tt  \|_2^2
 \right)   +
\frac{1}{R_2} \left(  \|\pp_z \eta \|_2^2  + \|\pp_z \tt \|_2^2
\right)    \right]  \; ds \leq K_7,   \label{K-7}
\end{eqnarray}
where
\begin{eqnarray}
&&\hskip-.68in K_7= C e^{(K_1+K_3^{2/3} +K_5^{2/3} +K_6^2) t} \left[
\|v_0\|_{H^1(\Om)}^2 + K_1 +  \|Q\|_2^2+ K_2^2\; (K_2+K_3^{1/3}
+K_5^{1/3} +K_6) \right].   \label{K7}
\end{eqnarray}

\subsection{$\|{\Dd}_{H} \overline{v}\|_{H^1(M)}^2 +
\|{\nabla}_{H}({\nabla}_H^{\perp} \cdot v_z)\|_{H^1(\Om)}^2 + \|
{\nabla}_{H}\left({\nabla}_{H} \cdot v_z + R_1 T\right)
\|_{H^1(\Om)}^2 + \|{\nabla}_{H} T\|_{H^1(\Om)}^2$ {\bf estimates}}

First, let us observe that
\begin{eqnarray*}
&&\hskip-.68in
|{\nabla}_{H} v(x,y,z)| \leq  |{\nabla}_{H}
\overline{v}(x,y)| +  \int_{-h}^0 |{\nabla}_{H} v_z(x,y,z)| \; dz.
\end{eqnarray*}
Therefore, from the above and (\ref{ZZZ}), we have
\begin{eqnarray*}
&&\hskip-.68in \|{\nabla}_{H} v\|_{\infty} \leq  \|{\nabla}_{H}
\overline{v}\|_{\infty} + \left\|\int_{-h}^0 |{\nabla}_{H} u| \; dz
\right\|_{\infty}  \\
&&\hskip-.68in
 \leq  \|{\nabla}_{H} \overline{v}\|_{\infty} + R_1 \int_{-h}^0 \|{\nabla}_{H} \beta\|_{\infty} \;
dz + \left\|\int_{-h}^0 |{\nabla}_{H} \zeta| \; dz
\right\|_{\infty}.
\end{eqnarray*}
By applying inequality  (\ref{BW-1}) to
${\nabla}_{H} \overline{v}$ and $\int_{-h}^0
|{\nabla}_{H} \zeta| \; dz$ we reach
\begin{eqnarray}
&&\hskip-.68in
\|{\nabla}_{H} \overline{v}\|_{\infty}
 \leq  C \|{\nabla}_{H} \overline{v}\|_{H^1(M)} \left(1+\log^+  \|{\nabla}_{H}
 \overline{v}\|_{H^2(M)}\right)^{1/2},  \label{EE-1} \\
&&\hskip-.58in
\left\|\int_{-h}^0 |{\nabla}_{H} \zeta| \; dz
\right\|_{\infty} \leq C \left\|\int_{-h}^0 |{\nabla}_{H} \zeta|
\; dz \right\|_{H^1(M)}\left(1+ \log^+ \left\|\int_{-h}^0
|{\nabla}_{H} \zeta| \; dz \right\|_{H^2(M)}\right)^{1/2}
\label{EE-2}  \\
\end{eqnarray}
Applying inequality (\ref{BW-2}) to ${\nabla}_{H}
\beta$, also by (\ref{BETA-1}) and (\ref{BETA-2}),
we reach
\begin{eqnarray}
&&\hskip-.68in  \|{\nabla}_{H} \beta\|_{\infty}
 \leq C \left( \|{\nabla}_{H} \cdot \beta \|_{\infty} + \|{\nabla}_H^{\perp} \cdot \beta \|_{\infty}  \right)
 \left(1+\log^+  \|{\nabla}_{H} \beta\|_{H^2(M)}\right)\leq C
\|T \|_{\infty} \left(1+\log^+ \|T\|_{H^2(M)}\right)   \label{EE-3}
\end{eqnarray}
Therefore, by (\ref{EE-1})--(\ref{EE-3}), we infer
that
\begin{eqnarray}
&&\hskip-.68in \|{\nabla}_{H} v\|_{\infty}
 \leq C \|{\nabla}_{H} \overline{v}\|_{H^1(M)} \left(1+\log^+  \|{\nabla}_{H}
 \overline{v}\|_{H^2(M)}\right)^{1/2} + C\int_{-h}^0 \left[
\|T \|_{\infty} \left(1+\log^+ \|T\|_{H^2(M)}\right) \right]\; dz
\nonumber \\
&&\hskip-.58in
 + C \left\|\int_{-h}^0 |{\nabla}_{H} \zeta|
\; dz \right\|_{H^1(M)}\left(1+ \log^+ \left\|\int_{-h}^0
|{\nabla}_{H} \zeta| \; dz \right\|_{H^2(M)}\right)^{1/2}
\nonumber \\
&&\hskip-.68in
 \leq  C \|{\nabla}_{H} \overline{v}\|_{H^1(M)} \left(1+\log^+  \|{\nabla}_{H}
 \overline{v}\|_{H^2(M)}\right)^{1/2}  + C \|T
\|_{\infty}  \left(1+\log^+ \int_{-h}^0 \|T\|_{H^2(M)} \; dz \right)
\nonumber \\
&&\hskip-.58in
 + C \left\|\int_{-h}^0 |{\nabla}_{H} \zeta|
\; dz \right\|_{H^1(M)}\left(1+ \log^+ \left\|\int_{-h}^0
|{\nabla}_{H} \zeta| \; dz \right\|_{H^2(M)}\right)^{1/2}
\nonumber  \\
&&\hskip-.68in
 \leq C \|{\nabla}_{H} \overline{v}\|_{H^1(M)} \left(1+\log^+  \|{\nabla}_{H}
 \overline{v}\|_{H^2(M)}\right)^{1/2}  + C \left\|T\right\|_{\infty} \left(1+\log^+
 \|{\Dd}_{H} T\|_{L^2(\Om)}\right)
 \nonumber \\
&&\hskip-.58in
 + C \left( \| \eta \|_{H^1(\Om)}+\| \tt \|_{H^1(\Om)}\right)
 \left[1+ \log^+ \left( \| {\Dd}_{H} \eta \|_{L^2(\Om)}+\|{\Dd}_{H} \tt
 \|_{L^2(\Om)}\right)\right]^{1/2}.   \label{LOGG}
\end{eqnarray}

\vskip0.1in

\subsubsection{$\|{\nabla}_{H} {\Dd}_{H} \overline{v}\|_2^2 $ {\bf estimates}}
By taking the ${\Dd}_{H}$ to  equation (\ref{EQ1})
and then taking the inner product of equation
(\ref{EQ1}) with $ - {\Dd}_{H}^2 \overline{v}$ in
$L^2(M)$, we reach
\begin{eqnarray*}
&&\hskip-.68in \frac{1}{2} \frac{d \|{\nabla}_{H} {\Dd}_{H}
\overline{v}
\|_{2}^{2} }{d t} + \frac{1}{R_1} \|{\Dd}_{H}^2 \overline{v} \|_2^2  \\
&&\hskip-.68in
 =  \int_{M} {\Dd}_{H} \left\{
(\overline{v} \cdot {\nabla}_{H} ) \overline{v} + \overline{ \left[
(\widetilde{v} \cdot {\nabla}_{H}) \widetilde{v}  + ({\nabla}_{H}
\cdot \widetilde{v}) \; \widetilde{v}\right]} +f_0 \vec{k}\times
\overline{v} \right\} \cdot
{\Dd}_{H}^2 \overline{v} \; dxdy.
\end{eqnarray*}
Integrating by parts and applying (\ref{EQ2}), we
obtain
\begin{eqnarray*}
&&\hskip-.68in \frac{1}{2} \frac{d \|{\nabla}_{H} {\Dd}_{H}
\overline{v}
\|_{2}^{2} }{d t} + \frac{1}{R_1} \|{\Dd}_{H}^2 \overline{v} \|_2^2  \\
&&\hskip-.68in \leq C \int_{M} \left\{|\overline{v}|\;
|{\nabla}_{H}^3 \overline{v}|+ |{\nabla}_{H} \overline{v}|\;
|{\nabla}_{H}^2 \overline{v}| + \int_{-h}^0 \left(
|\widetilde{v}|\;|{\nabla}_{H}^3 \widetilde{v}|+ |{\nabla}_{H}
\widetilde{v}|\; |{\nabla}_{H}^2 \widetilde{v}| \right)\; dz +
|{\Dd}_{H} \overline{v} | \right\}
|{\Dd}_{H}^2 \overline{v}| \; dxdy    \\
&&\hskip-.68in \leq C \int_{M} \left\{|\overline{v}|\;
|{\nabla}_{H}^3 \overline{v}|+ |{\nabla}_{H} \overline{v}|\;
|{\nabla}_{H}^2 \overline{v}| + \left( \int_{-h}^0
|\widetilde{v}|\;dz \right) \left( \int_{-h}^0 (|{\nabla}_{H}^3
\zeta| + |{\nabla}_{H}^3 \beta| ) \; dz\right)
\right.   \\
&&\hskip-.58in \left. +  \left( \int_{-h}^0 (|{\nabla}_{H}
\zeta|+|{\nabla}_{H} \beta|) \; dz \right) \left(  \int_{-h}^0 (
|{\nabla}_{H}^2 \zeta|+|{\nabla}_{H}^2 \beta| ) \; dz  \right) +
|{\Dd}_{H} \overline{v} | \right\} |{\Dd}_{H}^2 \overline{v}| \;
dxdy.
\end{eqnarray*}
By applying H\"{o}lder inequality, (\ref{SI-1}), (\ref{SI-111}), \ref{DIV-CUR},
(\ref{MAIN-1}) and (\ref{INQ}) to the above estimate, we obtain
\begin{eqnarray*}
&&\hskip-.168in \frac{1}{2} \frac{d \|{\nabla}_{H} {\Dd}_{H}
\overline{v}
\|_{2}^{2} }{d t} + \frac{1}{R_1} \|{\Dd}_{H}^2 \overline{v} \|_2^2    \\
&&\hskip-.268in \leq C  \left\{\|\overline{v}\|_4\; \|{\nabla}_{H}^3
\overline{v}\|_4+ \|{\nabla}_{H} \overline{v}\|_4\; \|{\nabla}_{H}^2
\overline{v}\|_4 + \int_{-h}^0\|\widetilde{v}\|_{\infty}\; dz
\left\|\int_{-h}^0
(|{\nabla}_{H}^3 \zeta| + |{\nabla}_{H}^3 \beta| ) \; dz \right\|_2 \right.   \\
&&\hskip-.058in \left. + \left\| \int_{-h}^0 (|{\nabla}_{H}
\zeta|+|{\nabla}_{H} \beta|) \; dz \right\|_4 \left\| \int_{-h}^0 (
|{\nabla}_{H}^2 \zeta|+|{\nabla}_{H}^2 \beta| ) \; dz  \right\|_4 +
\|{\Dd}_{H} \overline{v}\|_2 \right\}
\|{\Dd}_{H}^2 \overline{v}\|_2  \\
&&\hskip-.268in \leq C
\left\{\|\overline{v}\|_2^{\frac{1}{2}}\|{\nabla}_{H}
\overline{v}\|_2^{\frac{1}{2}}\; \|{\nabla}_{H}^3
\overline{v}\|_2^{\frac{1}{2}} \|{\Dd}_{H}^2
\overline{v}\|_2^{\frac{1}{2}} + \|{\nabla}_{H} \overline{v}\|_2\;
\|{\nabla}_{H}^3 \overline{v}\|_2 \right.   \\
&&\hskip-.058in \left.  + \|\widetilde{v}\|_2^{\frac{1}{2}}
\|{\nabla}_{H} \widetilde{v}\|_{\infty}^{\frac{1}{2}} \left(
\|{\Dd}_{H} \eta\|_2
+\|{\Dd}_{H} \tt\|_2 + \|{\Dd}_{H} T\|_2 \right) \right.   \\
&&\hskip-.058in \left. + \left( \|\eta\|_2+ \|\tt\|_2+\|T\|_{2}
\right) \left(  \|{\nabla}_{H}^2 \eta\|_2+\|{\nabla}_{H}^2
\tt\|_2+\|{\nabla}_{H}^2 T\|_2 \right) +\|{\Dd}_{H} \tt\|_2 \;\;
\right\} \|{\Dd}_{H}^2 \overline{v}\|_2.
\end{eqnarray*}
Thus, by Young's and the Cauchy--Schwarz inequalities, we
have
\begin{eqnarray}
&&\hskip-.68in \frac{d \|{\nabla}_{H} {\Dd}_{H} \overline{v}
\|_{2}^{2} }{d t} + \frac{1}{R_1} \|{\Dd}_{H}^2 \overline{v} \|_2^2  \nonumber
\\
&&\hskip-.58in
\leq C \left( \|\overline{v}\|_2^2\|{\nabla}_{H} \overline{v}\|_2^2
+ \|{\nabla}_{H} \overline{v}\|_2^2 \right) \; \|{\nabla}_{H}
{\Dd}_{H} \overline{v}\|_2^2 \nonumber
\\
&&\hskip-.58in + \left( \|\widetilde{v}\|_2 \|{\nabla}_{H}
\widetilde{v}\|_{\infty}   +  \|\eta\|_2^2+ \|\tt\|_2^2+\|T\|_{2}^2
\right) \left(  \|{\Dd}_{H} \eta\|_2^2+\|{\Dd}_{H}
\tt\|_2^2+\|{\Dd}_{H} T\|_2^2 \right). \label{VB-EST}
\end{eqnarray}

\subsubsection{$\|{\Dd}_{H} T\|_2+\|{\nabla}_{H} T_z\|_2$ {\bf estimates}}
By applying the operator  ${\Dd}_{H}$ to  equation
(\ref{EQ5}), and then taking the inner product of
equation (\ref{EQ5}) with $ {\Dd}_{H} T+T_{zz}$ in
$L^2(\Om)$,  we get
\begin{eqnarray*}
&&\hskip-.168in \frac{1}{2} \frac{d (\|{\Dd}_{H}
T\|_2^2+\|{\nabla}_{H} T_z\|_2^2)}{dt} +
\frac{1}{R_3}\left(\|{\Dd}_{H} T_z\|_2^2 +\|{\nabla}_{H} T_{zz}\|_2^2 \right) \\
&&\hskip-.165in = - \int_{\Om} {\Dd}_{H} \left[ v \cdot {\nabla}_{H}
T - \left( \int_{-h}^z {\nabla}_{H} \cdot v(x,y, \xi,t) d\xi
\right) \left(\frac{\pp T}{\pp z}+\frac{1}{h}\right) -Q \right] \; {\Dd}_{H} T \; dxdydz   \\
&&\hskip-.05in - \int_{\Om} {\nabla}_{H} \left[ v \cdot {\nabla}_{H}
T_z - \left( \int_{-h}^z {\nabla}_{H} \cdot v(x,y, \xi,t) d\xi
\right) \frac{\pp^2 T}{\pp z^2}+ u \cdot {\nabla}_{H} T  \right.  \\
&&\hskip-.05in \left.  - \left(
{\nabla}_{H} \cdot v \right) \left(\frac{\pp T}{\pp
z}+\frac{1}{h}\right)
-Q_z \right] \cdot {\nabla}_{H} T_z \; dxdydz \\
&&\hskip-.165in = - \int_{\Om}  \left[ {\Dd}_{H} v \cdot
{\nabla}_{H} T +2 {\nabla}_{H} v \cdot {\nabla}_{H}^2 T - \left(
\int_{-h}^z {\Dd}_{H}({\nabla}_{H} \cdot v(x,y, \xi,t)) d\xi \right)
\left(\frac{\pp T}{\pp
z}+\frac{1}{h}\right)   \right.  \\
&&\hskip-.05in \left.  - 2 \left( \int_{-h}^z {\nabla}_{H}
({\nabla}_{H} \cdot v(x,y, \xi,t)) d\xi
\right) {\nabla}_{H} T_z  -{\Dd}_{H} Q \right] \; {\Dd}_{H} T \; dxdydz   \\
&&\hskip-.05in - \int_{\Om}  \left[ {\nabla}_{H} v \cdot
{\nabla}_{H} T_z - \left( \int_{-h}^z {\nabla}_{H}({\nabla}_{H}
\cdot v(x,y, \xi,t)) d\xi \right)
\frac{\pp^2 T}{\pp z^2}  \right.  \\
&&\hskip-.05in \left. + {\nabla}_{H} u \cdot {\nabla}_{H} T +u \cdot
{\nabla}_{H}^2 T - {\nabla}_{H} \left( {\nabla}_{H} \cdot v \right)
\left(\frac{\pp T}{\pp z}+\frac{1}{h}\right)-\left( {\nabla}_{H}
\cdot v \right) {\nabla}_{H} T_z
-{\nabla}_{H} Q_z \right] \cdot {\nabla}_{H} T_z \; dxdydz  \\
&&\hskip-.165in = - \int_{\Om}  \left[ {\Dd}_{H} v \cdot
{\nabla}_{H} T +2 {\nabla}_{H} v \cdot {\nabla}_{H}^2 T -
\frac{1}{h}  \left( \int_{-h}^z {\Dd}_{H} ({\nabla}_{H} \cdot v(x,y,
\xi,t)) d\xi \right) + {\Dd}_{H}({\nabla}_{H} \cdot v) T
    \right.  \\
&&\hskip-.05in \left.  + 2 {\nabla}_{H} ({\nabla}_{H} \cdot v) \cdot
{\nabla}_{H} T  -{\Dd}_{H} Q \right] \; {\Dd}_{H} T \;
dxdydz  \\
&&\hskip-.05in - \int_{\Om}  \left[ \int_{-h}^z
{\Dd}_{H}({\nabla}_{H} \cdot v(x,y, \xi,t)) d\xi  \; T + 2
\int_{-h}^z {\nabla}_{H} ({\nabla}_{H} \cdot v(x,y, \xi,t)) d\xi
\cdot {\nabla}_{H} T \right] \cdot {\Dd}_{H} T_z \; dxdydz
 \\
&&\hskip-.05in - \int_{\Om}  \left[ {\nabla}_{H} v \cdot
{\nabla}_{H} T_z + \left(  {\nabla}_{H}({\nabla}_{H} \cdot v
)\right)
\frac{\pp T}{\pp z} + {\nabla}_{H} u \cdot {\nabla}_{H} T +u \cdot
{\nabla}_{H}^2 T   \right.  \\
&&\hskip-.05in \left.  - \frac{1}{h} {\nabla}_{H} \left( {\nabla}_{H}
\cdot v \right) + {\nabla}_{H} \left( {\nabla}_{H} \cdot u \right) T
+ \left( {\nabla}_{H} \cdot u \right) {\nabla}_{H} T
-{\nabla}_{H} Q_z \right] \cdot {\nabla}_{H} T_z \; dxdydz \\
&&\hskip-.05in - \int_{\Om}  \left[  \left( \int_{-h}^z
{\nabla}_{H}({\nabla}_{H} \cdot v(x,y, \xi,t)) d\xi \right)
\frac{\pp T}{\pp z}+  {\nabla}_{H} \left( {\nabla}_{H} \cdot v
\right) T + \left( {\nabla}_{H} \cdot v \right) {\nabla}_{H} T
\right] \cdot {\nabla}_{H} T_{zz} \; dxdydz.
\end{eqnarray*}
Thus,
\begin{eqnarray*}
&&\hskip-.168in \frac{1}{2} \frac{d (\|{\Dd}_{H}
T\|_2^2+\|{\nabla}_{H} T_z\|_2^2)}{dt} +
\frac{1}{R_3}\left(\|{\Dd}_{H} T_z\|_2^2 +\|{\nabla}_{H} T_{zz}\|_2^2 \right)
\\
&&\hskip-.168in \leq C \int_{\Om} \left\{ |{\nabla}_{H} v| |
{\nabla}_{H}^2 T|+ |{\Dd}_{H} v|\, |{\nabla}_{H}
T|+\overline{|{\Dd}_{H} ({\nabla}_{H} \cdot v)|}
+|{\nabla}_{H}^2({\nabla}_{H} \cdot  v)|\; |T | + |{\nabla}_{H}
({\nabla}_{H} \cdot v)|\; |{\nabla}_{H} T | \right.  \\
&&\hskip-.004in \left. +|{\Dd}_{H} Q| \right\} \; |{\Dd}_{H} T| \;
dxdydz
\\
&&\hskip-.05in  + C \int_{\Om}\left\{  \left( \int_{-h}^0 |{\Dd}_{H}
({\nabla}_{H} \cdot v)|\; dz |T |+\int_{-h}^0
|{\nabla}_{H}({\nabla}_{H} \cdot v)|\; dz |{\nabla}_{H} T |
\right) |{\Dd}_{H} T_z|  \right\} \; dxdydz  \\
&&\hskip-.05in + C \int_{\Om} \left\{ |{\nabla}_{H} v| |
{\nabla}_{H} T_z|+ | {\nabla}_{H}({\nabla}_{H} \cdot v )|\, |T_z| +
| {\nabla}_{H} u |\,  |{\nabla}_{H} T|+ |u|\;|{\nabla}_{H}^2 T|
 \right. \\
&&\hskip-.005in \left. +|{\nabla}_{H} ({\nabla}_{H} \cdot v)| +|{\nabla}_{H} ( {\nabla}_{H}
\cdot u )|\,  |T|  +|{\nabla}_{H} Q_z| \right\} \; |{\nabla}_{H}
T_z| \; dxdydz
\\
&&\hskip-.05in  + C \int_{\Om} \left[ \left( \int_{-h}^0
|{\nabla}_{H} ({\nabla}_{H}\cdot v)|\; dz \right) |T_z | +
|{\nabla}_{H}({\nabla}_{H} \cdot v)| |T | +  |{\nabla}_{H} \cdot v|
|{\nabla}_{H} T | \right] |{\nabla}_{H} T_{zz}|   \; dxdydz.
\end{eqnarray*}
Thanks to  (\ref{ZZZ})--(\ref{TTT}), we obtain
\begin{eqnarray*}
&&\hskip-.168in \frac{1}{2} \frac{d (\|{\Dd}_{H}
T\|_2^2+\|{\nabla}_{H} T_z\|_2^2)}{dt} +
\frac{1}{R_3}\left(\|{\Dd}_{H} T_z\|_2^2 +\|{\nabla}_{H} T_{zz}\|_2^2 \right)   \\
&&\hskip-.168in \leq C \int_{\Om} \left\{ |{\nabla}_{H} v| |
{\nabla}_{H}^2 T|+ \int_{-h}^0 (|{\Dd}_{H} \zeta|+|{\nabla}_{H}
T|)\;dz \, |{\nabla}_{H} T|+\int_{-h}^0 (|{\Dd}_{H}
\tt|+|{\Dd}_{H} T|)\; dz  \right.  \\
&&\hskip-.05in \left. +\int_{-h}^0 (|{\nabla}_{H}^2 \tt|+|{\nabla}_{H}^2 T|)\, dz \; |T |
+\int_{-h}^0 (|{\nabla}_{H} \tt|
 +|{\nabla}_{H} T|)\;dz |{\nabla}_{H} T |  +|{\Dd}_{H} Q| \right\} \; |{\Dd}_{H} T| \;
dxdydz
\\
&&\hskip-.05in  + C \int_{\Om}\left\{  \left( \int_{-h}^0
(|{\nabla}_{H}^2 \tt|+|{\nabla}_{H}^2 T|) \; dz |T |+\int_{-h}^0
(|{\nabla}_{H} \tt|+|{\nabla}_{H} T|)\; dz |{\nabla}_{H} T |
\right) |{\Dd}_{H} T_z|  \right\} \; dxdydz  \\
&&\hskip-.05in + C \int_{\Om} \left\{ |{\nabla}_{H} v| |
{\nabla}_{H} T_z|+ \int_{-h}^0 (|{\nabla}_{H} \tt|+|{\nabla}_{H}
T|)\,dz\; |T_z| + | {\nabla}_{H} u |\, |{\nabla}_{H} T| \right.  \\
&&\hskip-.004in \left. + |u|\;|{\nabla}_{H}^2 T| +
\int_{-h}^0(|{\nabla}_{H} \tt|+|{\nabla}_{H} T|)\; dz
+(|{\nabla}_{H} \tt|+|{\nabla}_{H} T|)\,  |T| +|{\nabla}_{H} Q_z|
\right\} \; |{\nabla}_{H} T_z| \; dxdydz
\\
&&\hskip-.05in  + C \int_{\Om} \left[ \left( \int_{-h}^0
(|{\nabla}_{H} \tt|+|{\nabla}_{H} T|)\; dz \right) |T_z | +
(|{\nabla}_{H} \tt|+|{\nabla}_{H} T|) |T | \right.   \\
&&\hskip-.004in \left.
+  \int_{-h}^0
(|\tt|+|T|)\; dz |{\nabla}_{H} T | \right] |{\nabla}_{H} T_{zz}|   \;
dxdydz.
\end{eqnarray*}
Using H\"{o}lder inequality, and inequalities
(\ref{MAIN-1}), (\ref{MAIN-2}) and (\ref{INQ}), we
obtain
\begin{eqnarray*}
&&\hskip-.168in \frac{1}{2} \frac{d (\|{\Dd}_{H}
T\|_2^2+\|{\nabla}_{H} T_z\|_2^2)}{dt} +
\frac{1}{R_3}\left(\|{\Dd}_{H} T_z\|_2^2 +\|{\nabla}_{H}
T_{zz}\|_2^2 \right)
 \\
&&\hskip-.168in \leq C \|{\nabla}_{H} v\|_{\infty} \| {\Dd}_{H}
T\|_2^2 + C \|{\Dd}_{H} \zeta\|_2^{1/2}\|{\nabla}_{H} {\Dd}_{H}
\zeta\|_2^{1/2} \|T\|_{\infty}^{1/2} \|{\nabla}_{H}^2 T\|_2^{3/2} +
C \|T\|_{\infty}
\|{\nabla}_{H}^2 T\|_2^2 \\
&&\hskip-.05in + C \|{\nabla}_{H}^2 \tt\|_2\,  (1+\|T \|_{\infty})
\| {\Dd}_{H} T\|_2
 + C \|{\nabla}_{H} \tt\|_2^{1/2}\|{\nabla}_{H}^2  \tt\|_2^{1/2}
\|T\|_{\infty}^{1/2} \|{\nabla}_{H}^2 T\|_2^{3/2}  + C \|{\Dd}_{H} Q\|_2  \; \|{\Dd}_{H} T\|_2 \\
&&\hskip-.05in
 + C (\|{\nabla}_{H}^2 \tt\|_2+\|{\nabla}_{H}^2 T\|_2) \;  \|T \|_{\infty} \|{\Dd}_{H} T_z\|_2
+C \|{\nabla}_{H} \tt\|_2^{1/2}\|{\nabla}_{H}^2  \tt\|_2^{1/2}
\|T\|_{\infty}^{1/2} \|{\nabla}_{H}^2 T\|_2^{1/2}  \|{\Dd}_{H}
T_z\|_2
\\
&&\hskip-.05in + C \|{\nabla}_{H} v\|_{\infty} \| {\nabla}_{H}
T_z\|_2^2 +
C \|{\nabla}_{H} \tt\|_2^{1/2}\|{\nabla}_{H}^2  \tt\|_2^{1/2} \|T_z\|_2^{1/2}\|{\nabla}_{H}  T_z\|_2^{3/2}  \\
&&\hskip-.05in + C \|T\|_{\infty}^{1/2} \|{\nabla}_{H}^2 T\|_2^{1/2}
\|T_z\|_2^{1/2}\|{\nabla}_{H}  T_z\|_2^{3/2}
 + C \| {\nabla}_{H} u \|_3\,
\|{\nabla}_{H} T\|_6 \|{\nabla}_{H} T_z\|_2 +
C \|u\|_{\infty} \;\|{\nabla}_{H}^2 T\|_2  \|{\nabla}_{H}  T_z\|_2 \\
&&\hskip-.05in +C  (\|{\nabla}_{H} \tt\|_2+\|{\nabla}_{H} T\|_2 )(1+ \|T\|_{\infty})   \\
&&\hskip-.05in
 + C \left( \|{\nabla}_{H} \tt\|_2^{1/2}\|{\nabla}_{H}^2  \tt\|_2^{1/2}
+C  \|T\|_{\infty}^{1/2} \|{\nabla}_{H}^2 T\|_2^{1/2} \right)
\|T_z\|_2^{1/2}\|{\nabla}_{H}  T_z\|_2^{1/2}
\|{\nabla}_{H}  T_{zz}\|_2 \\
&&\hskip-.05in
 +\|{\nabla}_{H} Q_z\|_2  \;
\|{\nabla}_{H} T_z\|_2 + C (\|{\nabla}_{H} \tt\|_2+\|{\nabla}_{H}
T\|_2)
\|T\|_{\infty}  \|{\nabla}_{H}  T_{zz}\|_2  \\
&&\hskip-.05in +C  \|\tt\|_2^{1/2}\|{\nabla}_{H}  \tt\|_2^{1/2}
\|T\|_{\infty}^{1/2} \|{\nabla}_{H}^2 T\|_2^{1/2} \|{\nabla}_{H}  T_{zz}\|_2   \\
&&\hskip-.168in \leq C \|{\nabla}_{H} v\|_{\infty} (\| {\Dd}_{H}
T\|_2^2+\| {\nabla}_{H} T_z\|_2^2)
 + C \left( \|{\nabla}_{H} \eta\|_2\ +\|{\nabla}_{H} \tt \|_2\right)^{1/2}
 \left( \|{\Dd}_{H} \eta\|_2\ +\|{\Dd}_{H} \tt \|_2\right)^{1/2}
\|T\|_{\infty}^{1/2} \|{\nabla}_{H}^2 T\|_2^{3/2}  \\
&&\hskip-.05in +C  \left(\|{\Dd}_{H} \eta\|_2 +\|{\Dd}_{H} \tt\|_2
+\|{\Dd}_{H} T\|_2 \right) \left(1+\|T\|_{\infty}\right) \;
\|{\Dd}_{H} T\|_2
 +C \|{\Dd}_{H} Q\|_2 \; \|{\Dd}_{H} T\|_2  +C \|{\nabla}_{H} Q_z\|_2 \; \|{\nabla}_{H} T_z\|_2  \\
&&\hskip-.05in
 +\left( \|{\nabla}_{H} \eta\|_2^{1/2}\|{\nabla}_{H}^2
\eta\|_2^{\frac{1}{2}}+\|{\nabla}_{H}
\tt\|_2^{\frac{1}{2}}\|{\nabla}_{H}^2
\tt\|_2^{\frac{1}{2}}+\|T\|_{\infty}^{\frac{1}{2}} \|{\Dd}_{H}
T\|_2^{\frac{1}{2}} \right) \; \|T\|_{\infty}^{\frac{1}{2}}
\|{\Dd}_{H} T\|_2^{\frac{1}{2}} \left(\|{\Dd}_{H} T\|_2+\|{\Dd}_{H}
T_z\|_2\right)  \\
&&\hskip-.05in   +C \|u\|_6  \|{\nabla}_{H}
T_z\|_2^{\frac{1}{2}}\|{\nabla}_{H}^2 T_z\|_2^{\frac{1}{2}}
\|{\Dd}_{H} T \| +C \left[ \|{\nabla}_{H}
\tt\|_2^{\frac{1}{2}}\|{\nabla}_{H}^2 \tt\|_2^{\frac{1}{2}}
\|T_z\|_2^{\frac{1}{2}} \|{\nabla}_{H}
T_z\|_2^{\frac{1}{2}} \right. \\
&&\hskip-.05in   + \left. \left(
\|T\|_4+\|\tt\|_2^{1/2}\|{\nabla}_{H} \tt\|_2^{1/2}+ \|T_z\|_2^{1/2}
\|{\nabla}_{H} T_z\|_2^{1/2} \right) \; \|T\|_{\infty}^{1/2}
\|{\Dd}_{H} T\|_2^{1/2}  \right] \|{\nabla}_{H} T_{zz}\|_2.
\end{eqnarray*}
By the Cauchy--Schwarz and Young's inequalities, we reach
\begin{eqnarray*}
&&\hskip-.68in \frac{d (\|{\Dd}_{H} T\|_2^2+\|{\nabla}_{H}
T_z\|_2^2)}{dt} +
\frac{1}{R_3}\left(\|{\Dd}_{H} T_z\|_2^2 +\|{\nabla}_{H} T_{zz}\|_2^2 \right)    \nonumber  \\
&&\hskip-.68in \leq  C\|{\nabla}_{H} v\|_{\infty} \left( \|{\Dd}_{H}
T\|_2^2+\|{\nabla}_{H} T_z\|_2^2 \right) + C \left( \|{\nabla}_{H}
\eta\|_2\ +\|{\nabla}_{H} \tt \|_2\right)^{1/2}\left( \|{\Dd}_{H}
\eta\|_2\ +\|{\Dd}_{H} \tt \|_2\right)^{1/2}
\|T\|_{\infty}^{1/2} \|{\Dd}_{H} T\|_2^{3/2}  \\
&&\hskip-.58in +C  \left(\|{\Dd}_{H} \eta\|_2 +\|{\Dd}_{H} \tt\|_2
+\|{\Dd}_{H} T\|_2 \right) \|T\|_{\infty} \; \|{\Dd}_{H} T\|_2
 +C \|{\Dd}_{H} Q\|_2 \; \|{\Dd}_{H} T\|_2  +C \|{\nabla}_{H} Q_z\|_2 \; \|{\nabla}_{H} T_z\|_2  \\
&&\hskip-.58in
 +\left( \|{\nabla}_{H} \eta\|_2^{1/2}\|{\Dd}_{H}
\eta\|_2^{\frac{1}{2}}+\|{\nabla}_{H}
\tt\|_2^{\frac{1}{2}}\|{\Dd}_{H}
\tt\|_2^{\frac{1}{2}}+\|T\|_{\infty}^{\frac{1}{2}} \|{\Dd}_{H}
T\|_2^{\frac{1}{2}} \right) \; \|T\|_{\infty}^{\frac{1}{2}}
\|{\Dd}_{H}
T\|_2^{\frac{1}{2}} \|{\Dd}_{H} T\|_2  \\
&&\hskip-.58in
 +\left( \|{\nabla}_{H} \eta\|_2\|{\Dd}_{H}
\eta\|_2+\|{\nabla}_{H} \tt\|_2\|{\Dd}_{H} \tt\|_2+\|T\|_{\infty}
\|{\Dd}_{H} T\|_2 \right) \; \|T\|_{\infty} \|{\Dd}_{H} T\|_2 +C
\|u\|_6^3 \left(\|{\nabla}_{H} T_z\|_2^2+  \|{\Dd}_{H} T
\|^2\right) \\
&&\hskip-.58in    +C \left[ \|{\nabla}_{H} \tt\|_2\|{\Dd}_{H}
\tt\|_2 \|T_z\|_2 \|{\nabla}_{H} T_z\|_2 + \left(
\|T\|_4^2+\|\tt\|_2\|{\nabla}_{H} \tt\|_2+ \|T_z\|_2 \|{\nabla}_{H}
T_z\|_2 \right) \; \|T\|_{\infty} \|{\Dd}_{H} T\|_2 \right]
\\
&&\hskip-.68in  \leq   C \|{\Dd}_{H} Q\|_2^2   +C \|{\nabla}_{H}
Q_z\|_2^2+ C \|T\|_{\infty}^{4} + C\|{\nabla}_{H} v\|_{\infty}
\left( \|{\Dd}_{H}
T\|_2^2+\|{\nabla}_{H} T_z\|_2^2 \right)  \\
&&\hskip-.58in + C \left(1+ \|{\nabla}_{H} \eta\|_2^2
+\|{\nabla}_{H} \tt \|_2^2 + \|T\|_{\infty}^{2}  +\|u\|_6^3
+\|T_z\|_2^2 \right) \left( \|{\Dd}_{H} T\|_2^2+\|{\nabla}_{H}
T_z\|_2^2 +\|{\Dd}_{H} \eta\|_2^2+\|{\Dd}_{H} \tt\|_2^2 \right).
\end{eqnarray*}
Thus, we get
\begin{eqnarray}
&&\hskip-.68in \frac{d (\|{\Dd}_{H} T\|_2^2+\|{\nabla}_{H}
T_z\|_2^2)}{dt} + \frac{1}{R_3}\left(\|{\Dd}_{H} T_z\|_2^2
+\|{\nabla}_{H} T_{zz}\|_2^2 \right) \nonumber
\\
&&\hskip-.68in  \leq   C \|{\Dd}_{H} Q\|_2^2   +C \|{\nabla}_{H}
Q_z\|_2^2+ C \|T\|_{\infty}^{4} + C\|{\nabla}_{H} v\|_{\infty}
\left( \|{\Dd}_{H}
T\|_2^2+\|{\nabla}_{H} T_z\|_2^2 \right)  \nonumber  \\
&&\hskip-.58in + C \left(1+ \|{\nabla}_{H} \eta\|_2^2
+\|{\nabla}_{H} \tt \|_2^2 + \|T\|_{\infty}^{2}  +\|u\|_6^3
+\|T_z\|_2^2 \right) \left( \|{\Dd}_{H} T\|_2^2+\|{\nabla}_{H}
T_z\|_2^2 +\|{\Dd}_{H} \eta\|_2^2+\|{\Dd}_{H} \tt\|_2^2 \right).
\label{T-EST}
\end{eqnarray}

\subsubsection{$\|{\nabla}_{H}({\nabla}_H^{\perp} \cdot v_z)\|_{H^1(\Om)}^2 + \| {\nabla}_{H}\left({\nabla}_{H}
\cdot v_z + R_1 T\right) \|_{H^1(\Om)}^2$ {\bf estimates}}
By acting with ${\Dd}_{H}$ on  equation (\ref{ETA}) and
equation (\ref{TT}), then taking the  inner product of  equation
(\ref{ETA}) with $ {\Dd}_{H} \eta +\eta_{zz}$ in $L^2$ and
equation (\ref{TT}) with ${\Dd}_{H} \tt+\tt_{zz}$ in $L^2$,
respectively, we get
\begin{eqnarray*}
&&\hskip-.18in \frac{1}{2} \frac{d \left( \|{\Dd}_{H}
\eta\|_2^2+\|{\nabla}_{H} \eta_{z}\|_2^2+\|{\Dd}_{H}
\tt\|_2^2+\|{\nabla}_{H} \tt_{z}\|_2^2\right)}{dt}   \\
&&\hskip-.08in + \frac{1}{R_1}
\left( \| {\nabla}_{H} {\Dd}_{H} \eta \|_2^2 +\| {\nabla}_{H}
{\Dd}_{H} \eta_z \|_2^2 + \| {\nabla}_{H} {\Dd}_{H} \tt \|_2^2+ \|
{\Dd}_{H} \tt_z \|_2^2
\right)\\
&&\hskip-.08in + \frac{1}{R_2}\left( \|{\Dd}_{H} \eta_z \|_2^2
+\|{\nabla}_{H} \eta_{zz} \|_2^2 + \|{\Dd}_{H} \tt_z \|_2^2
+ \|{\nabla}_{H} \tt_{zz} \|_2^2 \right)    \\
&&\hskip-.18in = \int_{\Om} {\nabla}_{H} \left\{ {\nabla}_H^{\perp} \cdot
\left[ (v \cdot {\nabla}_{H}) u - \left( \int_{-h}^z {\nabla}_{H}
\cdot v (x,y, \xi,t) d\xi \right) \frac{\pp u}{\pp z}  \right. \right.  \\
&&\hskip-.08in \left. \left. +(u \cdot {\nabla}_{H} ) v -
({\nabla}_{H} \cdot v) u  \right] +f_0\left( R_1\, T - \tt
\right)\right\} \cdot {\nabla}_{H} ({\Dd}_{H} \eta+\eta_{zz}) \;
dxdydz   \\
&&\hskip-.08in  +\int_{\Om} {\nabla}_{H} \left\{{\nabla}_{H} \cdot
\left[ (v \cdot {\nabla}_{H}) u - \left( \int_{-h}^z {\nabla}_{H}
\cdot v (x,y, \xi,t) d\xi \right) \frac{\pp u}{\pp z}  \right. \right.  \\
&&\hskip-.0in \left. \left. +  (u \cdot {\nabla}_{H} ) v -
({\nabla}_{H} \cdot v) u  \right]\, +f_0\,\eta \right\} \cdot
{\nabla}_{H} ({\Dd}_{H} \tt+\tt_{zz})
\; dxdydz   \\
&&\hskip-.08in - \int_{\Om} \left\{ R_1\left[ {\nabla}_{H} v \cdot
{\nabla}_{H} T+ v \cdot {\nabla}_{H}^2 T - \left( \int_{-h}^z
{\nabla}_{H} ({\nabla}_{H}
 \cdot v) (x,y, \xi,t) d\xi \right) \left( \frac{\pp T}{\pp
 z}+\frac{1}{h}\right)   \right. \right.   \\
&&\hskip-.08in  \left. \left.
 - {\nabla}_{H} Q +\, \left(\frac{1}{R_3}-\frac{1}{R_2}\right) {\nabla}_{H} T_{zz} \right] \,\right\}
  \cdot {\nabla}_{H} ({\Dd}_{H} \tt+\tt_{zz}) \; dxdydz
   \\
&&\hskip-.08in - \int_{\Om}  \left[  \left( \int_{-h}^z {\nabla}_{H}
({\nabla}_{H}
 \cdot v) (x,y, \xi,t) d\xi \right)  {\nabla}_{H} T_z  \; {\Dd}_{H} \tt  +
 \left( \int_{-h}^z {\nabla}_{H}
 \cdot v (x,y, \xi,t) d\xi \right)  {\nabla}_{H} T_z  \cdot {\nabla}_{H} \tt_{z}) \right] \; dxdydz
  \\
&&\hskip-.18in \leq C \int_{\Om} \left\{  \left(
|u|\;|{\nabla}_{H}^3 v| +|{\nabla}_{H} u|\;|{\nabla}_{H}^2 v|
+|{\nabla}_{H}^2 u|\;|{\nabla}_{H} v| +|{\nabla}_{H}^3 u|
|v| +|u_z|\;\int_{-h}^0|{\nabla}_{H}^2 ({\nabla}_{H} \cdot v)|\; dz   \right. \right.  \\
&&\hskip-.08in \left. \left.+|{\nabla}_{H}
u_z|\;\int_{-h}^0|{\nabla}_{H} ({\nabla}_{H} \cdot v)|\; dz
+|{\nabla}_{H}^2 u_z|\;\int_{-h}^0|{\nabla}_{H} \cdot v|\; dz
\right) \; \left(|{\nabla}_{H} {\Dd}_{H} \eta|+|{\nabla}_{H}
\eta_{zz}| +|{\nabla}_{H} {\Dd}_{H} \tt| +|{\nabla}_{H} \tt_{zz}|
\right) \right.   \\
&&\hskip-.08in \left. + C \left[ |{\nabla}_{H} v| |{\nabla}_{H} T| +
|v| |{\nabla}_{H}^2 T|+(1+|T_z|)\;\int_{-h}^0|{\nabla}_{H}
({\nabla}_{H} \cdot v)|\; dz
+|{\nabla}_{H} T_z|\;\int_{-h}^0|{\nabla}_{H} \cdot v|\; dz   \right. \right.   \\
&&\hskip-.08in  \left. \left. +|{\nabla}_{H} Q| +
\left|\frac{R_1}{R_3}-\frac{R_1}{R_2}\right| |{\nabla}_{H} T_{zz}|
\right] (|{\nabla}_{H} {\Dd}_{H} \tt|+|{\nabla}_{H} \tt_{zz}|)
\right\} \;
dxdydz  +C\left( \|{\Dd}_{H} T\|_2^2+\, \|{\Dd}_{H} \eta\|_2^2+\, \|{\Dd}_{H} \tt\|_2^2 \right) \\
&&\hskip-.18in \leq C \int_{\Om} \left\{ \left[ |u|\;\int_{-h}^0
(|{\nabla}_{H}^3 \zeta|+ |{\nabla}_{H}^3 \beta|) \; dz
+(|{\nabla}_{H} \zeta|+|{\nabla}_{H} \beta|) \;\int_{-h}^0
(|{\nabla}_{H}^2 \zeta|+ |{\nabla}_{H}^2 \beta|) \; dz +
(|{\nabla}_{H}^3 \zeta|+|{\nabla}_{H}^3 \beta|)|v| \right. \right.  \\
&&\hskip-.08in  +(|{\nabla}_{H}^2 \zeta|+|{\nabla}_{H}^2 \beta|)
\;\int_{-h}^0 (|{\nabla}_{H} \zeta|+ |{\nabla}_{H} \beta|) \; dz +
|u_z| \; \int_{-h}^0
(|{\nabla}_{H}^2 \tt|+ |{\nabla}_{H}^2 T|) \; dz  \\
&&\hskip-.08in +(|{\nabla}_{H} \zeta_z|+|{\nabla}_{H} \beta_z|)
\;\int_{-h}^0(|{\nabla}_{H} \tt|+|{\nabla}_{H} T|)\; dz  \\
&&\hskip-.08in \left. +(|{\nabla}_{H}^2 \zeta_z|+|{\nabla}_{H}^2
\beta_z|)\;\int_{-h}^0(|\tt|+|T|)\; dz \right] \;
\left(|{\nabla}_{H} {\Dd}_{H} \eta|+|{\nabla}_{H}  \eta_{zz}|
+|{\nabla}_{H} {\Dd}_{H} \tt| +|{\nabla}_{H} \tt_{zz}|
\right)    \\
&&\hskip-.08in  + C \left[ |{\nabla}_{H} T| \;
\int_{-h}^0(|{\nabla}_{H}
\zeta|+|{\nabla}_{H} \beta|)\; dz  + |v| |{\nabla}_{H}^2 T|  \right. \\
&&\hskip-.08in   +(1+|T_z|)\;\int_{-h}^0(|{\nabla}_{H}
\tt|+|{\nabla}_{H} T|)\; dz +|{\nabla}_{H}
T_z|\;\int_{-h}^0(|\tt|+|T|)\; dz    \\
&&\hskip-.08in  \left. \left. +|{\nabla}_{H} Q| +
\left|\frac{R_1}{R_3}-\frac{R_1}{R_2}\right| \; |{\nabla}_{H}
T_{zz}| \right] (|{\nabla}_{H} {\Dd}_{H} \tt|+|{\nabla}_{H}
\tt_{zz}|) \right\}\; dxdydz +C\left( \|{\Dd}_{H} T\|_2^2+\,
\|{\Dd}_{H} \eta\|_2^2+\, \|{\Dd}_{H} \tt\|_2^2 \right).
\end{eqnarray*}
Using H\"{o}lder inequality, and inequalities
(\ref{MAIN-1}), (\ref{MAIN-2}) and (\ref{INQ}), we
obtain
\begin{eqnarray*}
&&\hskip-.18in \frac{1}{2} \frac{d \left( \|{\Dd}_{H}
\eta\|_2^2+\|{\nabla}_{H} \eta_{z}\|_2^2+\|{\Dd}_{H}
\tt\|_2^2+\|{\nabla}_{H} \tt_{z}\|_2^2\right)}{dt}   \\
&&\hskip-.08in  + \frac{1}{R_1}
\left( \| {\nabla}_{H} {\Dd}_{H} \eta \|_2^2 +\| {\nabla}_{H}
{\Dd}_{H} \eta_z \|_2^2 + \| {\nabla}_{H} {\Dd}_{H} \tt \|_2^2+ \|
{\Dd}_{H} \tt_z \|_2^2
\right)  \\
&&\hskip-.08in  + \frac{1}{R_2}\left( \|{\Dd}_{H} \eta_z \|_2^2
+\|{\nabla}_{H} \eta_{zz} \|_2^2 + \|{\Dd}_{H} \tt_z \|_2^2 +
\|{\nabla}_{H} \tt_{zz} \|_2^2 \right)
\\
&&\hskip-.18in \leq C \left(\|{\nabla}_{H} {\Dd}_{H}
\eta\|_2+\|{\nabla}_{H} \eta_{zz}\|_2 +\|{\nabla}_{H} {\Dd}_{H}
\tt\|_2 +\|{\nabla}_{H} \tt_{zz}\|_2 \right)
 \left[  \left(\|u\|_{\infty} +\|v\|_{\infty}\right)\;
\left(\|{\nabla}_{H}^2 \eta\|_2+\|{\nabla}_{H}^2
\tt\|_2+\|{\nabla}_{H}^2 T\|_2\right)  \right.
\\
&&\hskip-.08in + \left(\| \eta\|_2^{\frac{1}{2}}\|{\nabla}_{H}
\eta\|_2^{\frac{1}{2}}+ \| \tt\|_2^{\frac{1}{2}}\|{\nabla}_{H}
\tt\|_2^{\frac{1}{2}}+ \|T\|_{\infty} \right) \;\left(\|{\nabla}_{H}
\eta\|_2^{\frac{1}{2}}\|{\nabla}_{H}^2 \eta\|_2^{\frac{1}{2}}+
\|{\nabla}_{H} \tt\|_2^{\frac{1}{2}}\|{\nabla}_{H}^2
\tt\|_2^{\frac{1}{2}}+ \| T\|_{\infty}^{\frac{1}{2}}\|{\nabla}_{H}^2
T\|_2^{\frac{1}{2}}\right)
\\
&&\hskip-.08in  +\left(\|{\nabla}_{H}
\eta\|_2^{\frac{1}{2}}\|{\nabla}_{H}^2 \eta\|_2^{\frac{1}{2}}+
\|{\nabla}_{H} \tt\|_2^{\frac{1}{2}}\|{\nabla}_{H}^2
\tt\|_2^{\frac{1}{2}}+ \|T\|_{\infty}^{\frac{1}{2}}\|{\nabla}_{H}^2
T\|_2^{\frac{1}{2}} \right) \left(\|
\eta\|_2^{\frac{1}{2}}\|{\nabla}_{H} \eta\|_2^{\frac{1}{2}}+ \|
\tt\|_2^{\frac{1}{2}}\|{\nabla}_{H} \tt\|_2^{\frac{1}{2}}+
\|T\|_{\infty} \right)
\\
&&\hskip-.08in
 + \left(\|\tt_z\|_2^{\frac{1}{2}}\|{\nabla}_{H} \tt_z\|_2^{\frac{1}{2}}
 + \|T_z\|_2^{\frac{1}{2}}\|{\nabla}_{H} T_z\|_2^{\frac{1}{2}} \right)
 \left(\|{\nabla}_{H} \tt\|_2^{\frac{1}{2}}\|{\nabla}_{H}^2 \tt\|_2^{\frac{1}{2}}
 + \|T\|_{\infty}^{\frac{1}{2}}\|{\nabla}_{H}^2 T\|_2^{\frac{1}{2}} \right)
 \\
&&\hskip-.08in  +\|u_z\|_2^{\frac{1}{2}}\|{\nabla}_{H}
u_z\|_2^{\frac{1}{2}} \;
 \left(\|{\nabla}_{H} \tt\|_2^{\frac{1}{2}}\|{\nabla}_{H}^2 \tt\|_2^{\frac{1}{2}}
 + \|T\|_{\infty}^{\frac{1}{2}}\|{\nabla}_{H}^2 T\|_2^{\frac{1}{2}} \right)
    \\
&&\hskip-.08in \left. +
 \left(\|{\nabla}_{H} \tt_z\|_2^{\frac{1}{2}}\|{\nabla}_{H}^2 \tt_z\|_2^{\frac{1}{2}}
 + \|{\nabla}_{H} T_z\|_2^{\frac{1}{2}}\|{\nabla}_{H}^2 T_z\|_2^{\frac{1}{2}} \right)
 \left(\|\tt\|_2^{\frac{1}{2}}\|{\nabla}_{H} \tt\|_2^{\frac{1}{2}}
 + \|T\|_{\infty} \right) \right] \;
   \\
&&\hskip-.08in + C \left[ \left(\|
\eta\|_2^{\frac{1}{2}}\|{\nabla}_{H} \eta\|_2^{\frac{1}{2}}+ \|
\tt\|_2^{\frac{1}{2}}\|{\nabla}_{H} \tt\|_2^{\frac{1}{2}}+ \|T\|_4
\right) \;\|T\|_{\infty}^{\frac{1}{2}}\|{\Dd}_{H}
T\|_2^{\frac{1}{2}} + \|v\|_{\infty} \|{\Dd}_{H} T\|_2
\right. \\
&&\hskip-.08in \left. +\left( \|{\nabla}_{H}
\tt\|_2^{\frac{1}{2}}\|{\nabla}_{H}^2 \tt\|_2^{\frac{1}{2}}+
\|T\|_{\infty}^{\frac{1}{2}} \|{\Dd}_{H} T\|_2^{\frac{1}{2}} \right)
\;\|T_z\|_2^{\frac{1}{2}}\|{\nabla}_{H} T_z\|_2^{\frac{1}{2}}
\right. \\
&&\hskip-.08in \left. +\left( \|\tt\|_2^{\frac{1}{2}}\|{\nabla}_{H}
\tt\|_2^{\frac{1}{2}}+ \|T\|_{\infty} \right) \;\|{\nabla}_{H}
T_z\|_2^{\frac{1}{2}}\|{\Dd}_{H}
T_z\|_2^{\frac{1}{2}}+\|{\nabla}_{H} \tt\|_2
+\|{\nabla}_{H} T\|_2 +\|{\nabla}_{H} Q\|_2 \right. \\
&&\hskip-.08in \left. + \left|\frac{R_1}{R_3}-\frac{R_1}{R_2}\right|
\|{\nabla}_{H} T_{zz}\|_2 \right] \left(\|{\nabla}_{H} {\Dd}_{H}
\tt\|_2+\|{\nabla}_{H} \tt_{zz}\|_2\right)+C\left( \|{\Dd}_{H}
T\|_2^2+\, \|{\Dd}_{H} \eta\|_2^2+\, \|{\Dd}_{H} \tt\|_2^2 \right)
\end{eqnarray*}
By Young's inequality and Cauchy--Schwarz inequality, we have
\begin{eqnarray}
&&\hskip-.148in  \frac{d \left( \|{\Dd}_{H}
\eta\|_2^2+\|{\nabla}_{H} \eta_{z}\|_2^2+\|{\Dd}_{H}
\tt\|_2^2+\|{\nabla}_{H} \tt_{z}\|_2^2\right)}{dt}  + \frac{1}{R_1}
\left( \| {\nabla}_{H} {\Dd}_{H} \eta \|_2^2 +\| {\nabla}_{H}
{\Dd}_{H} \eta_z \|_2^2 + \| {\nabla}_{H} {\Dd}_{H} \tt \|_2^2+ \|
{\Dd}_{H} \tt_z \|_2^2
\right) \nonumber  \\
&&\hskip-.038in  + \frac{1}{R_2}\left( \|{\Dd}_{H} \eta_z \|_2^2
+\|{\nabla}_{H} \eta_{zz} \|_2^2 + \|{\Dd}_{H} \tt_z \|_2^2
+ \|{\nabla}_{H} \tt_{zz} \|_2^2 \right)    \nonumber    \\
&&\hskip-.145in  \leq C \left(\|{\Dd}_{H} T\|_2^2+\|{\nabla}_{H}
T_{z}\|_2^2+ \|{\Dd}_{H} \eta\|_2^2+\|{\nabla}_{H}
\eta_{z}\|_2^2+\|{\Dd}_{H} \tt\|_2^2+\|{\nabla}_{H}
\tt_{z}\|_2^2\right) \left[1+ \|T\|_{\infty}^4+ \|{\nabla}_{H}
u_z\|^2_{2}+\|u_{zz}\|^2_{2}
\right. \nonumber  \\
&&\hskip-.038in \left. +\|{\nabla}_{H}
\eta\|^2_{2}+\|\eta_{z}\|^2_{2}
+\|{\nabla}_{H}\tt\|^2_{2}+\|\tt_{z}\|^2_{2}
 + \|T_z\|_2^2 +\|\eta\|_2^2 \|{\nabla}_{H} \eta\|_2^2 +\|\tt\|_2^2 \|{\nabla}_{H} \tt\|_2^2
\right]
\nonumber   \\
&&\hskip-.038in + \| \eta\|_2^2\|{\nabla}_{H} \eta\|_2^2+ \|
\tt\|_2^2 \|{\nabla}_{H} \tt\|_2^2 + \|T\|_{\infty}^4
 +\|u_z\|_2^2 \|{\nabla}_{H} u_z\|_2^2 + \|{\nabla}_{H} Q\|_2^2\nonumber   \\
&&\hskip-.038in + \frac{R_1^2(R_1+R_2)(R_2-R_3)^2}{R_2^2R_3^2}
\|{\nabla}_{H} T_{zz}\|_2^2.  \label{U-EST}
\end{eqnarray}
Next, we will obtain an estimate for
\[
\|{\nabla}_{H} {\Dd}_{H} \overline{v}\|_2^2 +
\|{\nabla}_{H}({\nabla}_H^{\perp} \cdot v_z)\|_{H^1(\Om)}^2 + \|
{\nabla}_{H}\left({\nabla}_{H} \cdot v_z + R_1 T \right)
\|_{H^1(\Om)}^2 + C_R \|{\nabla}_{H} T\|_{H^1(\Om)}^2,
\]
where $C_R = \frac{2R_1^2(R_1+R_2)(R_2-R_3)^2}{R_2^2R_3}.$
Denote by
\begin{eqnarray*}
&&\hskip-.68in  \mathcal{X}=1+\|{\nabla}_{H} {\Dd}_{H}
\overline{v}\|_2^2+C_R \, \|{\Dd}_{H} T\|_2^2+ C_R \, \|{\nabla}_{H}
T_{z}\|_2^2+ \|{\Dd}_{H} \eta\|_2^2+\|{\nabla}_{H}
\eta_{z}\|_2^2+\|{\Dd}_{H} \tt\|_2^2+\|{\nabla}_{H} \tt_{z}\|_2^2,    \\
&&\hskip-.68in  \mathcal{Y}=\|{\Dd}_{H}^2
\overline{v}\|_2^2+\|{\Dd}_{H} T_z\|_2^2+\|{\nabla}_{H}
T_{zz}\|_2^2+ \|{\nabla}_{H}{\Dd}_{H} \eta\|_2^2+\|{\Dd}_{H}
\eta_{z}\|_2^2+\|{\nabla}_{H}{\Dd}_{H} \tt\|_2^2+\|{\Dd}_{H}
\tt_{z}\|_2^2.
\end{eqnarray*}
Thus, by (\ref{LOGG}), we get
\[
\|{\nabla}_{H} v\|_{\infty} \leq C \left(  \|{\nabla}_{H}
\overline{v}\|_{H^1(M)}  +  \left\|T\right\|_{\infty} +\| \eta
\|_{H^1(\Om)}+\| \tt \|_{H^1(\Om)}  \right) \; \log \mathcal{X}.
\]
By virtue of  (\ref{VB-EST}), (\ref{T-EST}),
(\ref{U-EST}), Young's inequality, Cauchy--Schwarz
inequality and the above, we obtain
\begin{eqnarray*}
&&\hskip-.68in  \frac{d \mathcal{X}}{dt} + C\mathcal{Y} \leq C
\|v\|_2\;\left( \|{\nabla}_{H} \overline{v}\|_{H^1(M)}  +
\left\|T\right\|_{\infty} +\| \eta \|_{H^1(\Om)}+\| \tt
\|_{H^1(\Om)}  \right) \; \mathcal{X} \;
\log \mathcal{X}    \\
&&\hskip-.5in  + \left[1+ \|T\|_{\infty}^4+ \|T_z\|_2^2+
\|\overline{v}\|_2^2\left(1+\|\overline{v}\|_{H^1(M)}^2\right)   \right.    \\
&&\hskip-.5in   \left. +\|u_z\|^2_{H^1}
  +\left(1+\|\eta\|_2^2\right) \|\eta\|_{H^1}^2 +\left(1+\|\tt\|_2^2\right)
  \|\tt\|_{H^1}^2
\right] \; \mathcal{X}
   \\
&&\hskip-.5in + C \left\{ \|{\Dd}_{H} Q\|_2^2+\|{\nabla}_{H}
Q_z\|_2^2 + \| \eta\|_2^2\|{\nabla}_{H} \eta\|_2^2+ \| \tt\|_2^2
\|{\nabla}_{H} \tt\|_2^2 + \|T\|_{\infty}^4
 +\|u_z\|_2^2 \|{\nabla}_{H} u_z\|_2^2 \right\}.
\end{eqnarray*}
Let $\mathcal{X}=e^{\mathcal{Z}}.$ Then, $\frac{d \mathcal{X}}{dt}
= e^{\mathcal{Z}} \; \frac{d \mathcal{Z}}{dt}.$ As a result we
have
\begin{eqnarray}
&&\hskip-.68in  \frac{d \mathcal{Z}}{dt}  \leq C \|v\|_2\;\left(
\|{\nabla}_{H} \overline{v}\|_{H^1(M)}  +  \left\|T\right\|_{\infty}
+\|
\eta \|_{H^1(\Om)}+\| \tt \|_{H^1(\Om)}  \right) \; \mathcal{Z}  \nonumber   \\
&&\hskip-.5in  + \left[1+\|{\Dd}_{H} Q\|_2^2+\|{\nabla}_{H}
Q_z\|_2^2 + \|T\|_{\infty}^4+ \|T_z\|_2^2+
\|\overline{v}\|_2^2\left(1+\|\overline{v}\|_{H^1(M)}^2\right)   \right. \nonumber   \\
&&\hskip-.5in   \left. +\left(1+\|u_z\|_2^2\right) \|u_z\|^2_{H^1}
  +\left(1+\|\eta\|_2^2\right) \|\eta\|_{H^1}^2 +\left(1+\|\tt\|_2^2\right) \|\tt\|_{H^1}^2
\right].  \label{F-EST}
\end{eqnarray}
Thanks to Gronwall inequality, and the estimates
established in the previous subsections,  we get
\begin{eqnarray}
&&\hskip-.68in \mathcal{Z} \leq K,   \qquad  \mathcal{X} \leq e^K
   \label{K-F}
\end{eqnarray}
where
\begin{eqnarray}
&&\hskip-.68in K= e^{C(K_1+K_2+K_7)} \left[1+ \|v_0\|_{H^4}^2+\|
T_0\|_{H^2}^2 + t+\|{\Dd}_{H} Q\|_2^2\; t +\|{\nabla}_{H}
Q_z\|_2^2\; t \right]. \label{KF}
\end{eqnarray}
Moreover, we have
\begin{eqnarray}
&&\hskip-.68in \int_0^t \mathcal{Y} \; ds \leq K_8,
   \label{K-8}
\end{eqnarray}
where
\begin{eqnarray}
&&\hskip-.68in K_8= e^{C(K_1+K_2+K_7)} \left[1+ \|v_0\|_{H^4}^2+\|
T_0\|_{H^2}^2 + t+\|{\Dd}_{H} Q\|_2^2\; t +\|{\nabla}_{H}
Q_z\|_2^2\; t \right]. \label{K8}
\end{eqnarray}
Thanks to (\ref{K-6}) and the above we conclude
that the quantities $ \int_0^t
\|v_{zz}(s)\|_{H^1(\Om)}^2 \; ds$,
$\int_0^t\|{\Dd}_{H} \nabla_H v_z(s)\|_{H^1(\Om)}^2
\; ds $, $\|{\Dd}_{H} v_z(t)\|_{H^1(\Om)},
\|\nabla_H T (t)\|_{H^1(\Om)},$    and $ \int_0^t
\|\nabla_H T_z (s)\|_{H^1(\Om)}^2 \; ds$ are all
bounded uniformly in time, $t$, over the interval
$[0,\mathcal{T}_{*})$. Therefore, the strong
solution $(v(t), T(t))$ exists globally in time.
\end{proof}

\section{ Uniqueness of the Solutions} \label{S-4}

In this section  we state and prove the global
existence and uniqueness of the strong solution of
 system (\ref{EQV})--(\ref{EQ9}).

\begin{theorem} \label{T-MAIN}
Suppose that $Q \in H^2(\Om)$. Then for every $v_0
\in H^4(\Om)$, $T_0\in H^2(\Om)$ and
$\mathcal{T}>0,$ there is a unique solution $(v,
p_s, T)$ of  system {\em (\ref{EQV})--(\ref{EQ9})}
with
\begin{eqnarray*}
&& {\Dd}_{H} v_z, \;\; \nabla_H T \in C([0,\mathcal{T}], H^1(\Om)), \\
&& v_{zz}, \;\;{\Dd}_{H} \nabla_H v_z, \;\; {\nabla}_{H} T_z  \in L^2 ([0,\mathcal{T}], H^1(\Om)),
\end{eqnarray*}

\end{theorem}

\vskip0.05in

\begin{proof} In Theorem \ref{T-MAIN} of the previous section
we proved that the strong solutions exist globally
in time. Therefore, it remains to prove the
uniqueness of strong solutions, and their
continuous dependence on initial data, in the sense
specified by equation (\ref{LAST}) below. Let
$(v_1, (p_s)_1, T_1)$ and $(v_2, (p_s)_2, T_2)$ be
two strong solutions of system
(\ref{EQV})--(\ref{EQ9}) with initial data
$((v_0)_1, (T_0)_1)$ and $((v_0)_2, (T_0)_2)$,
respectively. Denote by $\phi=v_1-v_2, q_s =(p_s)_1
-(p_s)_2, \psi= T_1-T_2.$ It is clear that
\begin{eqnarray}
&&\hskip-.68in \frac{\pp \phi}{\pp t} + L_1 \phi + (v_1 \cdot
{\nabla}_{H}) \phi + (\phi \cdot {\nabla}_{H}) v_2 - \left(
\int_{-h}^z {\nabla}_{H} \cdot v_1(x,y, \xi,t) d\xi \right)
\frac{\pp \phi}{\pp  z} - \left( \int_{-h}^z
{\nabla}_{H} \cdot \phi(x,y, \xi,t) d\xi \right) \frac{\pp v_2}{\pp  z}  \nonumber  \\
&&\hskip-.6in +f_0 \vec{k}\times \phi + {\nabla}_{H} q_s -
{\nabla}_{H} \left( \int_{-h}^z
\psi (x,y,\xi,t) d\xi   \right)  =0,   \label{UEQ1}    \\
&&\hskip-.68in
 \frac{\pp \psi}{\pp t} -\frac{1}{R_3} \, \psi_{zz} + v_1 \cdot {\nabla}_{H} \psi + \phi
 \cdot {\nabla}_{H} T_2   \nonumber  \\
&&\hskip-.6in - \left( \int_{-h}^z {\nabla}_{H} \cdot v_1(x,y,
\xi,t)
 d\xi
\right) \frac{\pp \psi }{\pp z} - \left( \int_{-h}^z {\nabla}_{H}
\cdot \phi(x,y, \xi,t)
 d\xi
\right) \left(\frac{\pp T_2}{\pp z} +\frac{1}{h}\right)  = 0,
\label{UEQ2}
\end{eqnarray}
with initial data
\begin{eqnarray}
&&\hskip-.68in \phi(x,y,z,0) = (v_0)_1-(v_0)_2,  \label{UEQ3}   \\
&&\hskip-.68in \psi (x,y,z,0) = (T_0)_1-(T_0)_2. \label{UEQ4}
\end{eqnarray}
Taking the inner product of equation (\ref{UEQ1})
with  $\phi$ in $L^2(\Om)$, and equation
(\ref{UEQ2}) with $\psi$ in $L^2(\Om)$, we get
\begin{eqnarray*}
&&\hskip-.268in \frac{1}{2} \frac{d \|\phi\|_2^2}{dt}
+ \frac{1}{R_1} \|{\nabla}_{H} \phi\|_2^2 + \frac{1}{R_2}\|\phi_z\|_2^2  \\
&&\hskip-.265in =   - \int_{\Om} \left[ (v_1 \cdot {\nabla}_{H})
\phi + (\phi \cdot {\nabla}_{H}) v_2 - \left( \int_{-h}^z
{\nabla}_{H} \cdot v_1(x,y, \xi,t) d\xi \right) \frac{\pp \phi}{\pp
z} - \left( \int_{-h}^z {\nabla}_{H} \cdot \phi(x,y, \xi,t) d\xi
\right) \frac{\pp v_2}{\pp  z} \right] \cdot \phi
\; dxdydz    \\
&&\hskip-.16in
 + \int_{\Om} \left[  f_0 \vec{k}\times \phi + {\nabla}_{H} q_s -  {\nabla}_{H}  \left( \int_{-h}^z
\psi (x,y,\xi,t) d\xi   \right) \right] \cdot \phi \; dxdydz,
\end{eqnarray*}
and
\begin{eqnarray*}
&&\hskip-.168in \frac{1}{2} \frac{d \|\psi\|_2^2}{dt} +
\frac{1}{R_3}\|\psi_z\|_2^2 =   - \int_{\Om} \left[ v_1
 \cdot {\nabla}_{H} \psi + \phi
 \cdot {\nabla}_{H} T_2  \right. \\
&&\hskip-.16in  \left. - \left( \int_{-h}^z {\nabla}_{H} \cdot
v_1(x,y, \xi,t)
 d\xi
\right) \frac{\pp \psi }{\pp z} - \left( \int_{-h}^z {\nabla}_{H}
\cdot \phi(x,y, \xi,t)
 d\xi
\right) \left(\frac{\pp T_2}{\pp z} +\frac{1}{h}\right) \right]\; \psi \; dxdydz.
\end{eqnarray*}
Notice that
\begin{eqnarray}
&&\hskip-.65in  f_0 \vec{k}\times \phi   \cdot \psi  =0.
\label{DUT-1}
\end{eqnarray}
Integrating by parts, and using the boundary
conditions (\ref{EQ6}) and (\ref{EQ7}), we have
\begin{eqnarray}
&&\hskip-.65in   - \int_{\Om} \left( (v_1 \cdot {\nabla}_{H}) \phi -
\left( \int_{-h}^z {\nabla}_{H} \cdot v_1(x,y, \xi,t) d\xi \right)
\frac{\pp \phi}{\pp  z} \right) \cdot \phi \; dxdydz =0, \label{DUU-1}   \\
&&\hskip-.65in   - \int_{\Om} \left( v_1 \cdot {\nabla}_{H} \psi
 - \left( \int_{-h}^z {\nabla}_{H} \cdot v_1(x,y, \xi,t) d\xi \right)
\frac{\pp \psi}{\pp  z} \right) \cdot \psi \; dxdydz =0.
\label{DUT-11}
\end{eqnarray}
Integrating by parts, and using the boundary
conditions (\ref{EQ6}) and (\ref{EQ7}), we get
\begin{eqnarray}
&&\hskip-.65in   \int_{\Om} \left[   {\nabla}_{H} q_s - {\nabla}_{H}
\left( \int_{-h}^z \psi
(x,y,\xi,t) d\xi   \right) \right] \cdot \phi \; dxdydz \nonumber  \\
&&\hskip-.65in = \int_{\Om}   \left( \int_{-h}^z \psi (x,y,\xi,t)
d\xi \right) \; \left({\nabla}_{H} \cdot \phi \right) \; dxdydz.
\label{DUU-2}
\end{eqnarray}
Thus, by  (\ref{DUT-1})--(\ref{DUU-2}) we have
\begin{eqnarray*}
&&\hskip-.68in \frac{1}{2} \frac{d \|\phi\|_2^2}{dt} + \frac{1}{R_1}
\|{\nabla}_{H} \phi\|_2^2 + \frac{1}{R_2}\|\phi_z\|_2^2 = -
\int_{\Om} (\phi \cdot {\nabla}_{H}) v_2  \cdot \phi \;
dxdydz \\
&&\hskip-.65in  + \int_{\Om}  \int_{-h}^z {\nabla}_{H} \cdot
\phi(x,y, \xi,t) d\xi \frac{\pp v_2}{\pp  z} \cdot \phi \; dxdydz +
\int_{\Om} \left( \int_{-h}^z \psi (x,y,\xi,t) d\xi   \right) \;
\left({\nabla}_{H} \cdot \phi \right) \; dxdydz.
\end{eqnarray*}
and
\begin{eqnarray*}
&&\hskip-.68in \frac{1}{2} \frac{d \|\psi\|_2^2}{dt}+
\frac{1}{R_3}\|\psi_z\|_2^2 =   - \int_{\Om} (\phi \cdot
{\nabla}_{H}
T_2 ) \psi \; dxdydz \\
&&\hskip-.65in -\int_{\Om}  \left(  {\nabla}_{H} \cdot \phi \right)
T_2 \psi \; dxdydz- \int_{\Om}  \int_{-h}^z {\nabla}_{H} \cdot
\phi(x,y, \xi,t) d\xi  T_2  \psi_z \; dxdydz.
\end{eqnarray*}
Notice that by H\"{o}lder inequality and (\ref{MAIN-1})
\begin{eqnarray}
&&\hskip-.68in \left| \int_{\Om} (\phi \cdot {\nabla}_{H}) v_2 \cdot
\phi \; dxdydz \right| \leq  \left| \int_{\Om} |v_2|\,|\phi| \,
|{\nabla}_{H}\phi| \; dxdydz \right| \leq \| v_2\|_{6}
\|\phi\|_2^{\frac{1}{2}} \|{\nabla}_{H} \phi\|_2^{3/2},   \label{U1}
\\
&&\hskip-.68in \left | \int_{\Om} \phi \cdot {\nabla}_{H} T_2 \;
\psi \; dxdydz \right| \leq  \|{\nabla}_{H} T_2\|_3 \|\psi \|_2
\|\phi\|_6\leq C \|{\nabla}_{H} T_2\|_{H^1} \|\psi \|_2
\|{\nabla}_{H} \phi\|_2; \label{U-1}
\\
&&\hskip-.68in \left| \int_{\Om}  \int_{-h}^z {\nabla}_{H} \cdot
\phi(x,y, \xi,t) d\xi \frac{\pp v_2}{\pp  z} \cdot \phi \; dxdydz
\right|
\nonumber \\
&&\hskip-.68in  \leq C \|{\nabla}_{H} \phi\|_2 \left\|\frac{\pp
v_2}{\pp  z} \right\|_6 \|\phi\|_3  \leq C  \left\|\frac{\pp
v_2}{\pp  z} \right\|_6 \|\phi\|_2^{\frac{1}{2}}\|{\nabla}_{H}
\phi\|_2^{3/2}; \label{U3}
\\
&&\hskip-.68in \left| \int_{\Om}  \int_{-h}^z {\nabla}_{H} \cdot
\phi(x,y, \xi,t) d\xi \frac{\pp T_2}{\pp  z} \psi \; dxdydz \right|
\leq \left| \int_{\Om}  \int_{-h}^0 |{\nabla}_{H}\phi| dz
\int_{-h}^0 \left|\frac{\pp T_2}{\pp  z} \psi\right| dz  \; dxdy
\right|
\nonumber\\
&&\hskip-.68in \leq C \left| \int_{\Om}  \int_{-h}^0
|{\nabla}_{H}\phi| dz \left(\int_{-h}^0 \left|\frac{\pp T_2}{\pp
z}\right|^2\, dz\right)^{1/2} \left(\int_{-h}^0 |\psi|^2\,
dz\right)^{1/2} \; dxdy \right|
\nonumber  \\
&&\hskip-.68in \leq C \left\| \int_{-h}^0 \left|\frac{\pp T_2}{\pp
z}\right|^2\, dz \right\|_{\infty}^{1/2} \|{\nabla}_{H} \phi\|_2
\|\psi\|_2\leq C \left\|\frac{\pp {\Dd}_{H} T_2}{\pp z}\right\|_{2}
\|{\nabla}_{H} \phi\|_2 \|\psi\|_2. \label{U-3}
\end{eqnarray}
Therefore, by estimates (\ref{U1})--(\ref{U-3}), we reach
\begin{eqnarray*}
&&\hskip-.268in \frac{1}{2} \frac{d
\left(\|\phi\|_2^2+\|\psi\|_2^2\right) }{dt} + \frac{1}{R_1}
\|{\nabla}_{H} \phi\|_2^2 + \frac{1}{R_2}\|\phi_z\|_2^2
 + \frac{1}{R_3}\|\psi_z\|_2^2 \\
&&\hskip-.265in \leq C \left( \|v_2\|_{6} \|\phi\|_2^{\frac{1}{2}}
\|{\nabla}_{H} \phi\|_2^{3/2}+\|{\nabla}_{H} T_2\|_{H^1} \|\psi \|_2
\|{\nabla}_{H} \phi\|_2 + \| \pp_z v_2\|_6 \|\phi\|_2^{\frac{1}{2}}
\|{\nabla}_{H} \phi\|_2^{3/2} \right) + C \left\|\frac{\pp {\Dd}_{H}
T_2}{\pp z}\right\|_{2} \|{\nabla}_{H} \phi\|_2 \|\psi\|_2.
\end{eqnarray*}
By Young's inequality and the Cauchy--Schwarz inequality, we get
\begin{eqnarray*}
&&\hskip-.68in \frac{d \|\phi\|_2^2 +\|\psi(t)\|_2^2}{dt} \leq C
\left( \|v_2\|_{6}^4+\|{\nabla}_{H} T_2\|_{H^1}^2+ \| \pp_z
v_2\|_6^{4} + \left\|\frac{\pp {\Dd}_{H} T_2}{\pp z}\right\|_{2}^2
\right) \left( \|\phi\|_2^2 +\|\psi\|_2^2 \right).
\end{eqnarray*}
Thanks to Gronwall inequality, we obtain
\begin{eqnarray*}
&&\hskip-.68in \|\phi(t)\|_2^2 +\|\psi(t)\|_2^2
\leq \left( \|\phi(t=0)\|_2^2 +\|\psi(t=0)\|_2^2 \right) \times  \\
&&\hskip-.6in \exp \left\{ C  \int_0^t \left( \|v_2
(s)\|_{6}^4+\|{\nabla}_{H} T_2(s)\|_{H^1}^2+ \| \pp_z v_2(s)\|_6^{4}
+ \left\|\frac{\pp {\Dd}_{H} T_2}{\pp z}\right\|_{2}^2 \right) \; ds
\right\}.
\end{eqnarray*}
As a result of (\ref{K-F}),  we have
\begin{eqnarray}
&&\hskip-.68in \|\phi(t)\|_2^2 +\|\psi(t)\|_2^2 \leq \left(
\|\phi(t=0)\|_2^2 +\|\psi(t=0)\|_2^2 \right) \, \exp\left\{ C  \left(
(K_3^{2/3}  +K  + K_4^{2/3})\; t + K_8  \right)  \right\}.  \label{LAST}
\end{eqnarray}
The above inequality proves the continuous dependence of the
solutions on the initial data, and in particular, when
$\phi(t=0)=0$ and $\psi(t=0)=0$, we have $\phi(t)=0$ and $\psi(t)=0,$ for all
$t\ge 0$. Therefore, the strong solution is unique.

\end{proof}

\noindent
\section*{Acknowledgments}
The authors are thankful to the kind hospitality of
the Institute for Mathematics and its Applications
(IMA), University of Minnesota, where part of this
work was completed. This work was supported in part
by the NSF grants no.~DMS-0709228, DMS-0708832 and
DMS-1009950. The work of C.C. was also supported in
part by the Florida International University
Foundation. E.S.T. also acknowledges the warm
hospitality of the Freie Universit\"{a}t - Berlin
and  the partial support of the Alexander von
Humboldt Stiftung/Foundation and the Minerva
Stiftung/Foundation.

\end{document}